\documentclass[11pt,a4paper,leqno]{article}
\usepackage{amsmath}
\usepackage{amssymb}
\usepackage{graphicx}
\advance\topmargin by -2.2cm
\newcommand{\bbr}{I\!\!R}

\newcommand{\bbn}{I\!\!N}

\newcommand{\call}{{\cal L}}

\newcommand{\XXX}{{\bar X}}

\newcommand{\barr}{\begin{array}}
\newcommand{\earr}{\end{array}}
\newcommand{\beqq}{\begin{equation}}
\newcommand{\eeqq}{\end{equation}}
\newcommand{\beao}{\begin{eqnarray*}}
\newcommand{\eeao}{\end{eqnarray*}\noindent}
\newcommand{\beam}{\begin{eqnarray}}
\newcommand{\eeam}{\end{eqnarray}\noindent}

\newcommand{\halmos}{\quad\hfill\mbox{$\Box$}}

\newcommand{\al}{\alpha}

\newcommand{\vep}{\varepsilon}

\newtheorem{theo}{Theorem}
\newtheorem{prop}{\indent Proposition}

\newtheorem{defin}{\indent Definition}
\newtheorem{cor}{\indent Corollary}

\newcommand{\wt}{\widetilde}

\newcommand{\ov}{\overline} 
\newcommand{\ul}{\underline}

\author{R. H\"opfner$^*$ \and E.~L\"ocherbach \and M. Thieullen\thanks{This work has been supported by the Agence Nationale de la Recherche through the project MANDy, Mathematical Analysis of Neuronal Dynamics, ANR-09-BLAN-0008-01. e-mail addresses: \tt{hoepfner@mathematik.uni-mainz.de},
     \tt{eva.loecherbach@u-cergy.fr} and
    \tt{michele.thieullen@upmc.fr}} \\
{\it Johannes Gutenberg-Universit\"at Mainz, Universit\'e de Cergy-Pontoise} \\
 {\it and Universit\'e Pierre et Marie Curie.}}
\begin{document}

\title{Strongly degenerate time inhomogeneous SDEs: densities and support properties. Application to a Hodgkin-Huxley system \\ with periodic input.}

\maketitle

\begin{abstract}
In this paper we study the existence of densities for strongly degenerate stochastic differential equations (SDEs) whose coefficients depend on time and are not globally Lipschitz. In these models neither local ellipticity nor the strong H\"ormander condition is satisfied. In this general setting we show that continuous transition densities indeed exist in all neighborhoods of points where the weak H\"ormander condition is satisfied. We also exhibit regions where these densities remain positive. We then apply these results to stochastic Hodgkin-Huxley models with periodic input as a first step towards the study of ergodicity properties of such systems in the sense of \cite{MT1}-\cite{MT2}.

\end{abstract}

\\


\section{Introduction}
This paper belongs to a series of three articles (see also \cite{HLT-2}, \cite{HLT-3}) in which we carry a statistical study based on the ergodic theorem for multidimensional strongly degenerate and time inhomogeneous random systems with a view towards inference in neuroscience. In order to achieve our goal we have to prove ergodicity of such systems and it is natural to prove that their solutions possess Lebesgue densities as well as to address properties of the support of their law. This is the topic of the present paper. We establish that densities exist, are continuous and strictly positive, at least on suitable parts of the state space. The coefficients of our stochastic differential equations (SDEs) depend on time and are not globally Lipschitz. The noise is one dimensional and present in the first and last component only. In order to prove existence of densities, the H\"ormander condition which provides a sufficient condition of existence has become classical. To the best of our knowledge the papers available in the literature all use the H\"ormander condition in the strong form or even in the form of local ellipticity. The systems we study in the present paper do not fit into this framework.

Our motivation is to describe with probabilistic tools the long time behavior of a neuron embedded in a network in order to be able to estimate either parameters of the model or the underlying network activity or characteristics of the spike trains generated by the neuron. The network activity is present via the synaptic stimulation the neuron receives through its dendritic tree. We describe the neuron by the Hodgkin-Huxley model which is very well known in physiology. This system is notoriously mathematically difficult and may exhibit a collection of different behaviors when submitted to a deterministic periodic input. The synaptic stimulation we consider is a random input carrying a deterministic and periodic signal. We are interested in ergodicity for the process composed of the neuron on the one hand and the input it receives on the other hand. This results in a  five dimensional (5D) time inhomogeneous random system driven by a one dimensional Brownian motion present in the first and last component only. This system belongs to the general class of SDEs mentioned above. Because of this original motivation we find it natural to introduce an intermediate  family of models that we call {\it SDEs with internal variables and random input} which lies in between the general class of SDEs and our specific 5D-stochastic Hodgkin-Huxley. This family includes all conductance based models with synaptic input relevant for modeling in neuroscience (the Hodgkin-Huxley system is a conductance-based model). It is also relevant in biology and physics since it describes a population of individuals of different types represented by the internal variables, coupled by a global variable. Each individual can occupy two states, active/inactive, open/closed.  The transition rates between these two states depend on the global variable only. The deterministic system on which the SDEs are built is obtained as the limit of a sequence of Piecewise Deterministic Markov Processes when the number of individuals goes to infinity (cf. \cite{CDMR}, \cite{RGGMK}, \cite{FGRC}, \cite{PTW-10}) in the sense of the Law of Large Numbers. When we consider SDEs with internal variables and random input we consider that the population is infinite, we neglect the intrinsic noise related to finite size effects and we focus on the external noise received from the environment. For instance, when we model a neuron, the individuals in the population are the ion channels, the global variable is the membrane potential (see  \cite{PTW-10}). For systems in this family we provide an explicit discussion of the H\"ormander condition that we later illustrate numerically in the last section devoted to the 5D-stochastic Hodgkin-Huxley model. Then we exhibit regions where densities, if they exist, remain positive. These regions are related to neighborhoods of equilibrium points of the underlying deterministic system. The particular structure of SDEs with internal variables and random input, namely some linearity in the internal equations, plays a key role. We also show that with positive probability, the solution of these systems can imitate any deterministic evolution resulting from an arbitrary input, on an arbitrary interval of time. This property is used in \cite{HLT-3} to prove that the 5D-stochastic Hodgkin-Huxley system spikes infinitely often almost surely. As a first step towards recurrence and ergodicity, we show in the present paper that the 5D-stochastic Hodgkin-Huxley model possesses Lebesgue densities. However, it is not straightforward to check whether this specific model satisfies the H\"ormander condition. Hence we decided to conduct a numerical study. This is one more example of the difference between a general theoretical study and the application to a specific situation. We show that the H\"ormander condition is satisfied at certain stable equilibrium points or along a specific stable periodic orbit. Therefore, depending on the starting point, the 5D-stochastic Hodgkin-Huxley model possesses strictly positive densities either in small neighborhoods of such equilibrium points or of a periodic orbit. We use this result in two companion papers \cite{HLT-2} and \cite{HLT-3}, where we address periodic ergodicity and prove limit theorems. To the best of our knowledge, no other probabilistic study has been presented in the literature before. There are some simulation studies (see. e.g. \cite{PPM-05} and \cite{YWWL-01}), but not much seems to be known mathematically.

We use Malliavin Calculus and the H\"ormander condition (cf. \cite{Nualart}) to establish the existence of densities. H\"ormander sufficient condition ensures that the diffusion in the random system is actually strong enough even if the noise is visible only on a restricted number of components. It is satisfied when the Lie algebra generated by the coefficients of the SDE has full dimension and can be found under two forms: the strong form involving the diffusion coefficients only, and the weak form possibly including the drift coefficient. To the best of our knowledge, all existing results for time inhomogeneous SDEs require at least the strong form of H\"ormander condition (see \cite{C-M} and the references therein) which fails to hold in our general class of SDEs. We can only hope for the weak H\"ormander condition to be satisfied. Moreover in general this condition will hold only locally. SDEs satisfying local H\"ormander condition with locally smooth coefficients have been considered recently in a time homogeneous setting (cf. \cite{Bally}, \cite{stefano}, \cite{Fournier}). In these works the local H\"ormander condition is ensured by a local ellipticity assumption (hence these papers deal with the strong form of this condition). However our framework is technically more difficult since time homogeneity fails and since we must work with the weak form of H\"ormander condition, which holds only locally. In this general setting we show that continuous transition densities indeed exist in all neighborhoods of points where the weak H\"ormander condition is satisfied. We also prove that these densities are lower semi continuous (lsc) w.r.t. the starting point even if our system does not enjoy the Feller property. In order to do so we extend a localization argument and estimates of the Fourier transform introduced in \cite{Bally}, \cite{stefano} and \cite{Fournier}. 

The present paper paper is organized as follows. In section \ref{models} we first present our general SDEs and assumptions, then we introduce SDEs with internal variables and random input. Section \ref{density general} is devoted to proving the existence of densities locally for time inhomogeneous systems with locally Lipschitz coefficients. In section \ref{ivri} we explicit the weak H\"ormander condition and address the positivity of densities for SDEs with internal variables and random input. The section \ref{sec:main} of the paper is devoted to the 5D-stochastic Hodgkin-Huxley model (with some reminders on the deterministic system). In the appendix (section \ref{app}) we provide complementary proofs.

\section{Our Models}\label{models}

We now describe the systems of SDEs that we consider in a general frame together with the assumptions under which we will work. Then we present the subclass of SDEs with internal variables and random input.

Given an integer $ m \geq 1$, we write $ x = (x_1, \ldots , x_m) $ for generic elements of $\bbr^m .$ Let $\sigma $ be a measurable function from $\bbr^m $ to $\bbr^m $ and $b$ a smooth function from $ [0, \infty [ \times \bbr^m $ to $\bbr^m . $ For all $ x \in \bbr^m,$ we consider the SDE 
\begin{equation}\label{eq:equationms}
X_{i,t}= x_i + \int_0^t b_i ( s, X_s) ds + \int_0^t \sigma_i (X_s) d W_s , \; t \geq 0 , \; i = 1, \ldots , m, 
\end{equation}
and assume throughout this paper that a unique strong solution exists (at least up to some lifetime).
Here, $W$ is a one-dimensional Brownian motion and $ \sigma$ is identified with an $m \otimes 1-$matrix. We write $P_x$ for the probability measure under which the solution $X = (X_t)_{t \geq 0 } $ of (\ref{eq:equationms}) starts at $x.$ Note that the time dependence is in the drift only. We assume that (\ref{eq:equationms}) satisfies the following assumptions.

{\bf (H1)}  There exists an increasing sequence of compacts $ \; K_n \subset K_{n+1}$ of the form $K_n = [a_n, b_n]= \prod_{i = 1}^m [a_{n,i} , b_{n,i} ]$ where $ a_n = ( a_{n,1} , \ldots , a_{n,m} )$, such that for any $x \in \bigcup_n K_n , $ the unique strong solution to \eqref{eq:equationms} starting from $x$ at time $0$ satisfies that  $T_n := \inf \{ t : X_t \notin K_n \}  \to \infty  \mbox{ almost surely as } n \to \infty .$

{\bf (H2)} The coefficients of  (\ref{eq:equationms}) are locally smooth. Namely we assume that for all $n$, $\, \sigma \in C_b^\infty (K_n , \bbr^m )$ and for every multi-index $\beta \in \{ 0, \ldots, m \}^l, l \geq 1, \quad  b(t,x) +  \partial_{\beta}  b (t,x)  \mbox{ is bounded on  }  [ 0, T ] \times K_n$ for all $T>0$. Here $\partial_\beta =\frac{ \partial^l}{ \partial x_{\beta_1} \ldots \partial x_{\beta_l} } $ and we identify $x_0 $ with $t.$  

Notice that as a consequence of assumption {\bf (H1)} we could choose as state space of the process $ ( X_t, t \geq 0 )$ the set $E := \bigcup_{n } K_n.$ We will do this in some parts of the paper, e.g.\ in Sections 4 and 5.

Following our original motivation, if we want to model a neuron embedded in a network from which it receives an input through its dendritic tree, and able to activate ion channels modeled by the internal variables $2, \ldots , m-1,$ we consider a subclass of (\ref{eq:equationms}) with the following particular structure.
\begin{eqnarray}\label{internal1}
dX_{1,t}&=& F (X_{l,t}, \ 1\leq l\leq m-1)dt +dX_{m,t}\, , \\
dX_{i,t}&=&[-a_i(X_{1,t})X_{i,t}+b_i(X_{1,t})]dt, \, \, \, i = 2, \cdot, \cdot, \cdot , m-1\, ,\nonumber\\
dX_{m,t} &=& b_m(t,X_{m,t}) dt \;+\sigma (X_{m,t})dW_t\, .\nonumber 
 \end{eqnarray}
Note that the last component $X_m$ follows an autonomous equation. $X_m$ is a random external input to the underlying deterministic system
\begin{eqnarray}\label{subinternal}
dz_{1,t}&=& F (z_{l,t}, \ 1\leq l\leq m-1) dt \, ,  \\
dz_{i,t}&=&[-a_i(z_{1,t})z_{i,t}+b_i(z_{1,t})]dt, \, \, \, i = 2, \cdot, \cdot, \cdot , m-1.\nonumber
 \end{eqnarray}
This system can be interpreted as the limit of a sequence of stochastic ones in the sense of the Law of Large Numbers or Fluid Limit (cf \cite{PTW-10}): consider a population of individuals of $m-1$ different types, with $N$ individuals of each type, each individual being in two states (active or inactive, open or closed). The individuals are coupled by a global variable $z_1$: the transition rates between the two states depend on $z_1$ only ($a_i(z_1)$ and $b_i(z_1)$ for individuals of type $i$). When $N\rightarrow \infty$, the variable $z_{i,t}$ gives the probability that an individual of type $i$ is active at time $t$. The detailed form of these functions in the Hodgkin-Huxley system (cf. \cite{HH-52}) will be provided in Section \ref{sec:main}. In this well-known system, three types of `agents' are considered which are responsible for opening or closing of K$^+$ and Na$^+$ ion channels. In this particular model, we have three equations for internal variables corresponding to $m=5.$ $z_1$ describes the  membrane potential of the neuron, which can be observed. The $z_i,\, i = 2, 3, 4$ are the gating variables associated to specific ion channels located in the membrane, that are not observed. One may consider models which include still more types of ion channels admitting their specific number of different types of `agents', hence the interest to consider models \eqref{internal1} with general $m.$ 

By the general assumptions associated with \eqref{eq:equationms}, the coefficients $ F : \bbr^{m-1} \to \bbr , $ $a_i , b_i : \bbr \to \bbr $ for $ 2 \le i \le m-1 $ and $b_m : [0, \infty ] \times \bbr \to \bbr $ are smooth. In what follows, we suppose that the coefficients of \eqref{internal1} are such that {\bf (H1)} and {\bf (H2)} are satisfied. If we assume moreover that $ 0 \le b_i (z_1) \le a_i (z_1) $ and $a_i ( z_1) > 0 $ for all $ i = 2, \ldots , m - 1 , $ for all $ z_1 \in \bbr , $ then 
\begin{equation}\label{eq:ninfty}
y_{i ,\infty}(z_1) := \frac{b_i (z_1)}{a_i ( z_1) } , \; z_1 \in \bbr , 
\end{equation}
are equilibrium points of the internal equations when we keep the first variable fixed at constant value $z_1.$ In particular, introducing 
\begin{equation}\label{eq:5}
 F_\infty (z_1) := F(z_1, y_{2, \infty } ( z_1) , \ldots , y_{m-1, \infty } (z_1) ) , \, z_1 \in \bbr , 
\end{equation}
any point $( z_1, y_{2, \infty } ( z_1) , \ldots , y_{m-1, \infty } (z_1) )$ such that $ F_\infty ( z_1) = 0 $ is an equilibrium point of the system \eqref{subinternal}. 


\section{Existence and smoothness of densities for (\ref{eq:equationms})}\label{density general}

\noindent Classically, one proves that the solution of an SDE admits a smooth density via Malliavin Calculus, imposing the H\"ormander condition. Many authors assume that the coefficients of the SDE are $C^\infty$, bounded, with bounded derivatives of any order and that the H\"ormander condition is satisfied all over the state space. Here the coefficients of (\ref{eq:equationms}) are not globally Lipschitz. Regarding the H\"ormander condition there are actually two possibilities: either to work under the strong H\"ormander condition or under the weak one which is a less stringent assumption. Many authors work under the strong H\"ormander condition. The strong degeneracy of (\ref{eq:equationms}) imposes to work under the weak form of H\"ormander condition, which moreover may hold only locally. In addition the drift coefficient depends on time. Hence we have to apply local arguments in a time inhomogeneous setting. 


\subsection{Local H\"ormander condition in a time dependent setting}\label{sec:hoerm} 
In this section we state our local weak H\"ormander condition. Localization has been considered in \cite{KS}, however in a time homogeneous framework. Later it was pointed out in \cite{C-M} that time dependence in the coefficients of an SDE may lead to serious difficulties. In what follows we describe the effect of time dependence on (\ref{eq:equationms}). The first step is to rewrite (\ref{eq:equationms})
in Stratonovich form. This amounts to replace the drift $ b(t,x) $ by 
\begin{equation}\label{driftstratonovich}
 \tilde b_i (t,x) := b_i (t,x) - \frac12 \sum_{k=1}^m \sigma_k (x) \frac{ \partial \sigma_i }{\partial x_k } (x) , \, \, \, \; 1 \le i \le m \; , \, \, \, x \in \bbr^m\, , 
 \end{equation}
which is still time inhomogeneous. Now let us extend the coefficients $\tilde b$ and $\sigma$ into vector fields (or linear differential operators of order one) on $[0,+\infty[\times {\bbr}^m$ by setting   
\begin{eqnarray*}
A_0 &=& \frac{\partial }{\partial t } + \sum_{ i = 1}^m \tilde b_i (t,x) \frac{\partial }{\partial x_i } = \frac{\partial}{\partial t} + \tilde b,\\
A_1&=& \sum_{ i =1}^m \sigma_i (x) \frac{\partial }{ \partial x_i }.
\end{eqnarray*}
$A_0$ and $A_1$ can be identified respectively with the $ (m+1 )-$dimensional function $ A_0 ( t, x ) = ( {\tilde b}_0, \tilde b_1 , \ldots , \tilde b_m) $ with ${\tilde b}_0=1$ and $ A_1 ( t, x ) = ( \sigma_0, \sigma_1 , \ldots , \sigma_m)$ with $\sigma_0=0$. Actually any vector field $ {\cal T}(t,x) ={\cal T}_0 (t,x) \frac{\partial }{\partial t } +\sum_{i=1}^m {\cal T}_i (t,x) \frac{\partial}{\partial x_i}$ can be identified with the $ (m+1 )-$dimensional function $ ( {\cal T}_0, {\cal T}_1 , \ldots , {\cal T}_m) $.

The Lie bracket of two vector fields $ {\cal T}(t,x) ={\cal T}_0 (t,x) \frac{\partial }{\partial t } +\sum_{i=1}^m {\cal T}_i (t,x) \frac{\partial}{\partial x_i}$ and $ {\cal V}(t,x) ={\cal V}_0 (t,x) \frac{\partial }{\partial t } +\sum_{i=1}^m {\cal V}_i (t,x) \frac{\partial}{\partial x_i}$ is defined as
\begin{equation*}
 [{\cal T}, {\cal V} ]_i := \sum_{j = 0}^m ( {\cal T}_j \frac{ \partial {\cal V}_i }{\partial x_j } - {\cal V}_j \frac{ \partial {\cal T}_i }{\partial x_j } )=({\cal T}_0 \frac{ \partial {\cal V}_i }{\partial t } - {\cal V}_0 \frac{ \partial {\cal T}_i }{\partial t})+\sum_{j =1}^m ( {\cal T}_j \frac{ \partial {\cal V}_i }{\partial x_j } - {\cal V}_j \frac{ \partial {\cal T}_i }{\partial x_j }) .
\end{equation*}
In particular since $\sigma$ in (\ref{eq:equationms}) does not depend on time,
\begin{equation*}
[A_1, {\cal V}]_i= \sum_{j =1}^m ( {\sigma}_j \frac{ \partial {\cal V}_i }{\partial x_j } - {\cal V}_j \frac{ \partial {\sigma }_i }{\partial x_j }),
\end{equation*}
where no time derivative appears. On the contrary a time derivative is present in $[A_0, {\cal V}]$ since
\begin{equation*}
 [{A_0}, {\cal V} ]_i=(\frac{ \partial {\cal V}_i }{\partial t } - {\cal V}_0 \frac{ \partial {\tilde b}_i }{\partial t})+\sum_{j =1}^m ( {\tilde b}_j \frac{ \partial {\cal V}_i }{\partial x_j } - {\cal V}_j \frac{ \partial {\tilde b}_i }{\partial x_j }) .
\end{equation*}

Notice that whenever $ {\cal V}_0$ vanishes identically, $[A_1, {\cal V}]_0\equiv0$ since $\sigma_0$ is zero and $[A_0, {\cal V}]_0\equiv 0$ as well since ${\tilde b}_0$ is constant equal to $1$. In this case the vector fields $[A_1, {\cal V}]$ and $[A_0, {\cal V}]$ belong to the $m$-dimensional space generated by the $\frac{\partial }{\partial x_i }, \, \, 1\leq i\leq m$. In particular $[A_1,A_0]$ belongs to this latter space as well as $A_1$ by definition.

Given $A_1$ and $A_0$ we can build two Lie algebras. On one hand the Lie algebra $\Lambda$ generated by the set $\{A_1,A_0\}$. On the other hand the Lie algebra $\call$ generated by the set $\{ A_1,[A_1,A_0]\}$. As just noticed the dimension of $\call(t,x)$ cannot exceed $m$ whatever $(t,x)\in [0,+\infty[\times {\bbr}^m$. However the dimension of $\Lambda (t,x)$ can be equal to $m+1$. Actually the following result holds. \begin{prop}\label{algebras}
For all $(t,x)\in [0,+\infty[\times {\bbr}^m, \, \, $ ${\rm dim}\, \,  \Lambda(t,x)={\rm dim}\, \,  \call(t,x)+1.$
\end{prop} 
Before giving the proof of this proposition we state the local H\"ormander condition we are going to work with and make an important remark about it. Recall $ E = \bigcup_{n} K_n \subset \bbr^m $ from condition {\bf (H1)}.

{\bf (LWH)} We say that the H\"ormander condition is satisfied at $(t,y_0)$ if there exist $r\in ]0,t[\,$ and $R>0$ such that $B_{5R}(y_0)\subset E$ and ${\rm dim}\, \, \Lambda(s,y)=m+1,\, \, \, \forall (s,y) \in [t-r,t]\times B_{3R}(y_0)$ (local weak H\"ormander condition).

A global version of {\bf (LWH)} has been introduced in \cite{C-M} as condition (1.6). This paper recalls an example from \cite{T} which points out the necessity to incorporate the operator $\frac{\partial}{\partial t}$ to the original framework via extension of the coefficient $\tilde b$ into $A_0$ described above. We refer the reader to \cite{C-M} for more details.

\noindent{\bf Proof of Proposition \ref{algebras}.} In practice when one wants to evaluate the dimension of an algebra like $\Lambda(t,x)$ or $\call(t,x)$, one computes iteratively Lie brackets of successive orders. Accordingly let us build the following sequences of Lie algebras indexed by $N\in\bbn$. Fix such an integer $N$. Construct first the set $\call_N$ iteratively such that it contains $A_1$ (initialization) and is stable by Lie brackets with $A_0$ and $A_1$ (iteration) of order up to $N$. Then define ${\rm LA}(\call_N)$ as the Lie algebra spanned by $\call_N$. Construct also iteratively the set $\Lambda_N$ such that it contains $A_1$ and $A_0$ (initialization) and is stable by Lie brackets with $A_0$ and $A_1$ (iteration) of order up to $N$. Then define ${\rm LA}(\Lambda_N)$ as the Lie algebra spanned by $\Lambda_N$. The difference in the initialization between $\call_N$ and $\Lambda_N$ plays a key role. It implies that $\Lambda_N\subset \{A_0\}\cup\call_N\cup {-\call_N}$ where we denote by $-\call_N$ the set $\{-L; L\in \call_N\}$. Moreover by definition ${\rm LA}(\Lambda_N)$ coincides with $\Lambda$ for all $N$. From these two facts we deduce that the dimension of ${\rm LA}(\call_N)$ does not depend on $N$ and that the identity holds.  \halmos 

In the sequel we will check {\bf (LWH)} at $(t,y_0)$ by successive computations of Lie brackets looking for $r\in ]0,t[\,$ and $N\in \bbn$ such that ${\rm dim}\, \, {\rm LA}(\call_N)(s,y_0)=m,\, \, \, \forall s \in [t-r,t]$ (as we did in the proof of Proposition \ref{algebras}).



  


\subsection{Local densities for (\ref{eq:equationms})}\label{localdensities}

Let us recall that an $\bbr^m $- valued random vector $Z$ admits a density with respect to Lebesgue measure or is absolutely continuous {\it on an open set $O \subset \bbr^m $}, if for some function $ p \in L^1 (O)$, 
\begin{equation*}
 E (f(X)) = \int f (y) p(y) dy ,
 \end{equation*}
 for any continuous and bounded function $f \in C_b ( \bbr^m ) $ satisfying $ supp (f) \subset O.$
\begin{theo}\label{theo:main2}
Assume that (\ref{eq:equationms}) satisfies {\bf (H1)} and {\bf (H2)}. Assume moreover that {\bf (LWH)} is satisfied at $(t,y_0).$ Fix $x \in {\bbr}^m $ and denote by $ ( X_t, t \geq 0 )$ the strong solution of (\ref{eq:equationms}) starting from $x$. Then the random variable $X_t $ admits a density on $ B_r ( y_0 ) $ which is continuous on $B_R ( y_0 ) , $ where $R$ is given in {\bf (LWH)}.
\end{theo}

Note that this density might be $ \equiv 0$ near $y_0;$ so far it is not granted that the process at time $t$ visits such neighborhoods for $ y_0 \in {\rm int} ( E) $ for arbitrary choice of a starting point $x \in \bbr^m $ with positive probability.

\begin{theo}\label{theo:lsc} Let us keep the assumptions and notations of Theorem \ref{theo:main2} and for $x$ in ${\bbr}^m$denote by $p_{0, t } (x,\cdot)$ the density of $X_t$ on $B_R ( y_0 )$. For any fixed $y \in B_R ( y_0) ,$ the map $x\in {\bbr}^m \mapsto  p_{0, t } (x,y) $ is lower semi-continuous.
\end{theo} 
Given the assumptions {\bf (H1)}-{\bf (H2)} on (\ref{eq:equationms}), we have to use localization arguments in order to prove these theorems. Localization arguments have been used in \cite{KS} and \cite{stefano}, however in a time homogeneous framework. Moreover in  \cite{stefano} the non degeneracy is local ellipticity which fails to hold for (\ref{eq:equationms}). We prove below that \cite{KS} and  \cite{stefano} can be extended to a time inhomogeneous SDE satisfying only the local weak H\"ormander condition {\bf (LWH)}.

\noindent{\bf Proof of Theorem \ref{theo:main2}.} In this proof we rely on the following criterion based on Fourier transform which ensures existence and regularity of Lebesgue densities. Let $\mu$ be a probability measure on $\bbr^m $ and $ \hat \mu$ its Fourier transform defined by $\hat \mu (\xi):= \frac{1}{(2 \pi)^m} \int_{\bbr^m} e^{ i < y,\xi >} \mu (y) dy$. If  $\hat \mu$
 is integrable, then $\mu$ is absolutely continuous and a continuous version of its density is given by 
\begin{equation}\label{inversefourier}
p(y) = \frac{1}{(2 \pi)^m} \int_{\bbr^m} e^{ - i < \xi , y>} \hat \mu (\xi ) d \xi. 
\end{equation}
Consider $R>0$ provided by {\bf (LWH)}. Denote by $\Phi $ a localizing function in $C_b^\infty ( \bbr^m ) $ satisfying $ 1_{ B_R (0)} \le \Phi \le 1_{ B_{2R}(0)}.$ Fix $x$, $t$ and $T$ with $t\leq T$ and assume that $ E_x ( \Phi ( X_t - y_0 )) := m_0 > 0 .$ We prove below that the probability measure $\nu$ defined by 
\begin{equation}\label{probanu}
 \int f(y) \nu (dy) : = \frac{1}{m_0} E_x \left( f( X_t) \Phi ( X_t - y_0 )\right), 
\end{equation}
for all $f \in C_b ( \bbr^m )$, is such that $\hat\nu$, given by $\hat \nu (\xi ) = \frac{1}{m_0} E_x \left( e^{i < \xi , X_t>} \Phi ( X_t - y_0 ) \right)$, is integrable. For this purpose the main step is to prove (\ref{eq:classical1})-(\ref{eq:classical2}) below. Although the form of (\ref{eq:classical1})-(\ref{eq:classical2}) is classical, we have to make sure that they hold in our time inhomogeneous framework. Let $\psi \in C^\infty_b (\bbr^m ) $ such that 
$$ \psi (y ) = \left\{ 
\begin{array}{ll} 
y & \mbox{ if } |y| \le 4 R \\
5 R \frac{ y}{|y|} & \mbox{ if } |y| \ge 5 R 
\end{array}
\right. $$
and $|\psi (y)| \le 5 R $ for all $ y .$ Let $ \bar b (t, y) = b (t, \psi ( y - y_0 ))$ and $\bar \sigma (y) = \sigma ( \psi ( y - y_0 ))$ be the localized coefficients of (\ref{eq:equationms}). Assumption {\bf (H2)} ensures that $\bar b$ and $\bar \sigma$ are $C_b^\infty-$extensions (w.r.t. $x$) of $ b_{| B_{4R} (y_0)} $ and $\sigma_{ | B_{4R} (y_0)} $ with $\bar b$ and its derivatives bounded on $[0,T]$. Let $ \XXX $ satisfy the SDE \begin{equation}\label{eq:processgood}
d\bar X_{i,s} = \bar b_i( s, \bar X_s) ds +  \bar \sigma_i (\bar X_s) d W_s ,\;  \, \, s \le T ,\; \, \, 1 \le i \le m,
\end{equation}
and $\bar X_{i,0}=X_0=x$. If $ x \in B_{4R } ( y_0), $ the processes $ \XXX $ and $ X$ coincide up to the first exit time of $ B_{4R} (y_0)$. For a fixed $ \delta$ in $] 0, t/2 \wedge r[, $ where $r$ is provided by {\bf (LWH)}, define $\tau_1 := \inf \{ s \geq t - \delta ;  X_s \in B_{3R} (y_0 ) \} $ and $\tau_2 := \inf \{ s \geq \tau_1 ; X_s \notin B_{ 4 R }(y_0) \} .$ 
The set $\{ \Phi ( X_t - y_0 ) > 0 \} $ is equal to the union 
\begin{equation*}
\{ \Phi ( X_t - y_0 ) > 0 ; t - \delta = \tau_1 < t < \tau_2  \} 
\cup \left\{ \Phi ( X_t - y_0 ) > 0 ; \sup_{ 0 \le s \le \delta } |  \bar X_{\tau_1 , \tau_1 + s} ( X_{\tau_1 } ) - X_{ \tau_1} | \geq R \right\},
\end{equation*}
where $ \bar X_{u,v} (z) $ denotes the value at time $v$ of the solution of (\ref{eq:processgood}) satisfying $\bar X_u=z$ at time $u$ when $u\leq v$ (classical notation for flows). Note that $\Phi(X_t-y_0)>0$ implies $X_t\in B_{2R}(y_0)$. Using the Markov property in $\tau_1$, we obtain the following expression of $\hat \nu$,
\begin{eqnarray*}\label{eq:fourier}
m_0 \hat \nu ( \xi ) &=&E_x \left( e^{ i < \xi, X_t>} \Phi ( X_t - y_0 )1_{\{ \Phi ( X_t - y_0 ) > 0\}}\, 1_{\{ \sup_{ 0 \le s \le \delta } | \XXX_{\tau_1, \tau_1 + s } ( X_{\tau_1}) - X_{ \tau_1} | \geq R\} }  \right)\\
&+& E_x \left( e^{ i < \xi, X_t>} \Phi ( X_t - y_0 )1_{\{  \Phi ( X_t - y_0 ) > 0\}}\, 1_{ \{ t- \delta = \tau_1 < t < \tau_2\}} \right).
\end{eqnarray*}
We are looking for upper bounds of $|\hat \nu ( \xi ) |$ to check whether $\hat \nu$ is integrable. The latter identity reads $m_0 \hat \nu ( \xi ) =A+B$. We will see shortly that the important contribution comes from $|B|$. To control $|A|$ we use the classical estimate
\begin{equation}\label{eq:control11}
P_x\left( \Phi ( X_t - y_0 ) > 0 ; \sup_{ 0 \le s \le \delta } | \XXX_{\tau_1, \tau_1 + s } ( X_{\tau_1}) - X_{ \tau_1} | \geq R \right) \le 
C(T,q,m , b, \sigma ) R^{-q} \delta^{q /2}.
\end{equation}
It is valid for all $ q > 0$ and holds uniformly in $x.$ The constant $C(T,q ,m,  b ,  \sigma ) $ depends  on the supremum norms of $\bar b $ and $\bar \sigma,$ hence by construction, on the supremum norms of $\sigma $ (resp. $b$) on $B_{5R} (y_0 )$ (resp. $B_{5R} (y_0 )\times [0,T]$). Notice that the right hand side of (\ref{eq:control11}) follows from\begin{equation}\label{eq:ub1}
E_x\left( \sup_{ u :s \le u \le t } | \bar X_{i,u} - \bar X_{i,s} |^q \right) \le C(T,q ,m,  b ,  \sigma ) (t-s)^{q/2}, \, \, \mbox{ for all $ 0 \le s \le t \le T.$}
\end{equation}
Let us now estimate $|B|$. Thanks to the Markov property at time $ t - \delta,$
\begin{equation}\label{eq:control2}
\Big| B \Big| \le 
 \sup_{ y \in B_{3R} (y_0)} | E_x \left( e^{ i < \xi, \XXX_{t-\delta, t} (y) > } \Phi ( \XXX_{t-\delta, t} (y)  - y_0 ) \right) | ,
\end{equation}
which again holds uniformly in $x.$ As in \cite{stefano}, we take advantage of the relationship between the exponential $ e^{ i < \xi, z >}$ and its second order partial derivatives w.r.t. each component of $z$. Namely $\partial^{(2)}_{z_\ell} e^{ i < \xi, z >}  =-\xi_\ell^2 e^{ i < \xi , z > }$. We denote by $\partial_\beta$ the composition of these second order partial derivatives for $\ell\in\{1,\cdot, \cdot,\cdot, m\}$ and set $ \| \xi \| := \prod_{\ell=1}^m  | \xi_\ell |$. Then
\begin{equation*}
| E_x \left( e^{ i < \xi , \XXX_{t-\delta, t} (y) >} \Phi ( \XXX_{t-\delta, t} (y)  - y_0) \right) | \leq 
\| \xi \|^{-2}  \Big| E_x \left( \partial_\beta e^{ i < \xi, \XXX_{t-\delta, t} (y)>} \Phi ( \XXX_{t-\delta, t} (y)- y_0 ) \right) \Big|.
\end{equation*}
From the integration by parts of Malliavin calculus we conclude that for some functional ${\cal H}$,
\begin{equation}\label{eq:classical1}
| E_x ( e^{ i < \xi , \XXX_{t-\delta, t} (y) >} \Phi ( \XXX_{t-\delta, t} (y)  - y_0) ) | \leq \| \xi \|^{-2} E_x( | {\cal H}( \XXX_{t-\delta, t} (y) , \Phi ( \XXX_{t-\delta, t} (y) - y_0 ) ) |). 
\end{equation}
We show in the Appendix (section \ref{app}) that  
\begin{equation}\label{eq:classical2}
\|  {\cal H} ( \XXX_{t-\delta, t} (y) , \Phi ( \XXX_{t-\delta, t} (y) - y_0 ))\|_p \le C (r,  p, R, m)  \delta^{ - m k_N  }.
\end{equation}
The constant $k_N$ depends on the order $N$ of successive Lie brackets needed to span $\bbr^m $ at any point of $B_{3R} (y_0 )$ according to {\bf (LWH)} (according to the remark after the proof of Proposition \ref{algebras}). We deduce from (\ref{eq:control11}) and (\ref{eq:classical2}) that, for any $ q \geq 1$ and any $0 <  \delta < \frac{t}{2} \wedge r ,$  
\begin{equation*}
m_0 |\hat \nu (\xi )| \le C (T, r, R, q, m ) \;
 \left[ R^{- q } \delta^{ q/2} + \| \xi \|^{-2} \delta^{ - m k_N} \right].
\end{equation*}
In order to bound $|\hat \nu |$ above by an integrable function we now exploit (as in \cite{stefano}) the freedom that still remains in the choice of the pair $(\delta, q)$. Indeed, for a given $\xi$, we can choose ($\delta, q)$ such that $ R^{- q } \delta^{ q/2} +\| \xi \|^{-2 }  \delta^{ - m k_N }$ tends to zero faster than $ \| \xi \|^{ - 3/2 }$ as $\| \xi \| \to \infty$ as follows:
\begin{equation*}
 \delta = t/2 \wedge r \wedge \| \xi \|^{-\frac{1}{2mk_N} }, \; \, \, \, q = 6 m  k_N. 
 \end{equation*}
Then $m_0 |\hat \nu (\xi ) |\le C (T, r, R, q, m ) \; \| \xi \|^{ - \frac32}$ and $\hat \nu$ is integrable. We conclude that (\ref{inversefourier}) holds for $\nu$. Therefore the solution of (\ref{eq:equationms}) starting from $x$ admits the density
\begin{equation}\label{densite}
p_{0, t } ( x, y )=\frac{1}{(2 \pi)^m} \int_{\bbr^m } e^{ - i <\xi , y > } E_x ( e^{ i < \xi , X_t>} \Phi ( X_t - y_0 ) ) d \xi,
\end{equation}
on $B_R (y_0)$. It remains to prove the continuity of $ p_{ 0, t} (x, y ) $ with respect to $y$. We split the integrals in (\ref{densite}) in two parts, over the bounded set $I:=\{ | \xi \| \le M \}$ and its complement $I^c $ for some $M>0$ and we apply the dominated convergence theorem. The modulus of the integrand is bounded on $I$. On $I^c$, we use the fact that $E_x ( e^{ i < \xi , X_t>} \Phi ( X_t - y_0 ) )$ coincides with $\hat\nu(\xi)$ and the inequality just established: $m_0 |\hat \nu (\xi ) |\le C (T, r, R, q, m ) \; \| \xi \|^{ - \frac32}.$ The continuity is uniform in $x$ since the upper bounds in (\ref{eq:control11}) and (\ref{eq:control2}) do not depend on $x.$\halmos

\noindent{\bf Proof of Theorem \ref{theo:lsc}.} We keep the notations introduced in the proof of Theorem \ref{theo:main2}, in particular $\Phi$ and $\nu$. In order to prove the lower semi-continuity w.r.t. $x$, it is enough to show that for fixed $y \in B_R (y_0)$, the function $p_{0, t} (\cdot,y)$ is the limit of an increasing sequence of continuous functions $x \mapsto p_{0, t}^{(n)} (x,y) $. We also use localization arguments here but now the approximating sequence is obtained by considering $X$ before it exits each compact $K_n$ (cf. {\bf (H1})). Note that continuous dependence on the starting point holds for each approximating process which enjoys the flow property whereas this property mail fail to hold for $X$ itself. So, given an integer $n$, let $b^{(n)} (t,x) $ and $\sigma^{(n)} (x)$ denote $C^\infty-$extensions (in $x$) of $ b(t, \cdot _{ | K_n} )$ and $\sigma_{| K_n} .$ Let $X^{(n)}$ be the solution of the localized version of (\ref{eq:equationms}) with coefficients $b^{(n)}$ and $\sigma^{(n)}$. The first exit time of $K_n$ by $X$ is denoted by $T_n$ (cf {\bf (H1)}). Using that $T_n \to \infty$ , we can write for any $x \in K_n$ and any positive measurable function $f,$ 
\begin{equation*}
m_0\int f(y) \nu (dy)=\lim_n \uparrow E_x \left( f( X_t)  \Phi (X_t - y_0 ) 1_{ \{ T_n > t \} } \right). 
\end{equation*}
Then for all $n$, since $X_t^{(n)} = X_t$ on $\{ T_n > t \} $ almost surely and $\Phi $ is non negative,
\begin{equation*}
m_0\int f(y) \nu (dy)\geq E_x \left( f( X_t ) \Phi (X_t - y_0 )  1_{ \{ T_n > t \}} \right) =E_x \left( f( X^{(n)}_t ) \Phi (X_t^{(n)} - y_0 )  1_{ \{ T_n > t \}} \right) .
\end{equation*}
We approximate $1_{ \{ T_n > t \}} $ by some continuous functional on $ \Omega := C( \bbr_+ , \bbr^m ) .$  The set $\Omega $ is endowed with the topology of uniform convergence on compacts. $\mathbb{P}_{0,x}^{(n)}$ denotes the law of $X^{(n)} $ on $ (\Omega , {\cal B} (\Omega ) ),$ starting from $x$ at time $0.$ The family $\{ \mathbb{P}_{0,x}^{(n)}, x \in \bbr^m \}$ has the Feller property, i.e. if $x_k \to x, $ then $\mathbb{P}_{0, x_k}^{(n)} \to \mathbb{P}_{0,x}^{(n)} $ weakly as $k \to \infty .$  Thanks to this property, $ E_x \left( f( X^{(n)}_t) \Phi (X_t^{(n)} - y_0 )   \right)$ is continuous w.r.t. $x$. Define $M_t^n = \max_{ s \le t } X_s^{(n)} $ and $m_t^n = \min_{s \le t} X_s^{(n)} $ coordinate-wise. Due to the structure of the compacts $K_n$ (see assumption {\bf (H1)}), we can construct  
$C^\infty-$functions $ \varphi^n, \Phi^n $ such that $ 1_{ [a_{n-1}, \infty [ } \le \varphi^n \le 1_{ [a_n,  \infty[ } $ and 
$ 1_{ ] - \infty , b_{n-1}] } \le \Phi^n \le 1_{ ] - \infty , b_n] } $ (these inequalities have to be understood coordinate-wise). Then, since $X_t$ equals $X_t^{(n)} $ up to time $T_n,$ 
\begin{equation*}
 \{ T_{n-1} > t \} = \{ a_{n -1} \le m_t^n \le M_t^n \le b_{ n - 1} \} \subset \{ \varphi^n( m_t^n ) > 0 , \Phi^n (M_t^n ) > 0  \} \subset \{ T_n > t \},
  \end{equation*}
and for any $ f \geq 0 , $ 
\begin{eqnarray*}
E_x \left(f( X^{(n)}_t)  \Phi (X_t^{(n)} - y_0 ) 1_{ \{ T_n > t \} } \right) 
& \geq & E_x \left( f( X^{(n)}_t) \Phi (X_t^{(n)} - y_0 )  \Phi^n ( M_t^n) \varphi^n (m_t^n) \right).
\end{eqnarray*}
Define now a sub-probability measure $\nu_n $ by 
\begin{equation}\label{eq:mun}
 m_0\int f(y) \nu_n ( dy ) := E_x \left( f( X^{(n)}_t ) \Phi (X_t^{(n)} - y_0 )  \Phi^n ( M_t^n) \varphi^n (m_t^n ) \right). 
\end{equation}
The new functional $ \Phi (X_t^{(n)} - y_0 )  \Phi^n ( M_t^n) \varphi^n (m_t^n )$ satisfies the same hypotheses as the former $ \Phi (X_t^{(n)} - y_0 ) . $
For any $f \geq 0 , $
\begin{equation*}
 \int f(y) \nu_n (dy) \le \int f(y) \nu_{n+1} (dy) \uparrow \int f(y) \nu (dy ) \mbox{ as } n \to \infty .
 \end{equation*} 
If we can show that $ \nu_n $ possesses a density, that we shall denote by $m_0^{ - 1 }  p_{0, t}^{(n)} (x, y ) ,$ the following inequalities will hold true  
\begin{equation}\label{eq:increase}
 p_{0, t}^{(n)} (x, y ) \le p_{0, t}^{(n+1)} (x, y ) \le p_{0, t} (x, y ) \quad \mbox{ for all $ n \geq 1 ,$}  
\end{equation}
for any fixed $ x ,$ $\lambda (dy) -\mbox{almost surely.}$
So in a next step we show that indeed $\nu_n $ possesses a density. In order to indicate explicitly the dependence on the starting point $x$, we introduce the notation $\gamma_n ( x, \xi)$ for $\hat \nu_n ( \xi )$ as follows,
\begin{equation*}
\gamma_n ( x, \xi) := \frac{1}{m_0}  E_x \left( e^{i < \xi, X^{(n)}_t>} \Phi (X_t^{(n)} - y_0 )  \Phi^n ( M_t^n) \varphi^n (m_t^n ) \right), 
\end{equation*}
and we apply the argument in the proof of Theorem \ref{theo:main2}. Inequalities (\ref{eq:control11}) and (\ref{eq:classical1})-(\ref{eq:classical2}) hold for $m_0 \gamma_n (x, \xi)$ which also satisfies $m_0 |\gamma_n (x, \xi) |\le C (T, r, R, q, m ) \; \| \xi \|^{ - \frac32}$. Therefore $\xi \to \gamma_n ( x, \xi ) $ is integrable. Hence, $ m_0 \nu_n $ admits a density that we denote $p_{0, t}^{(n)} (x, y ) $ given by
\begin{equation}\label{eq:pn}
p_{0, t}^{(n)} (x, y ) =  \frac{m_0}{(2\pi)^m} \int_{\bbr^m } e^{ - i < \xi , y >} \gamma_n ( x, \xi ) d \xi .
\end{equation}
From the fact that $ \gamma_n ( x, \xi ) \to \hat \nu ( \xi )$ as $ n \to \infty$ and that the upper bounds for $|\gamma_n|$ do not depend on $n$, we deduce that $p_{0, t}^{(n)} (x, y ) \to p_{0, t}  (x, y ).$ Taking into account  (\ref{eq:increase}), we conclude that $ p_{0, t} (x,y) = \lim_n \uparrow p_{0, t}^{(n)} (x,y) .$

It remains to show (by dominated convergence) that for any $y \in B_R ( y_0),$ the map $x \mapsto p_{0, t}^{(n)} (x,y) $ is continuous. This is a consequence of the continuity of $\gamma_n ( x, \xi) $ in $x $ (which follows from the Feller property of $\mathbb{P}_{0,x}^{(n)} $ and the fact that all operations appearing in $ \gamma_n (x, \xi) $ are continuous on $\Omega$) and the fact that (\ref{eq:control11}) and (\ref{eq:classical1})-(\ref{eq:classical2}) hold uniformly in $x.$ \halmos

\section{Densities for SDEs with internal variables and random input}\label{ivri}

In this section we consider the subclass of (\ref{eq:equationms}) presented in (\ref{internal1}), Section \ref{models}, that we call SDEs with internal variables and random input. For such systems we first make {\bf (LWH)} explicit. Then, when densities exist in the sense of the previous section, we address the question of their positivity. We show that with positive probability, the solution of (\ref{internal1}) can imitate any deterministic evolution resulting from an arbitrary input, on an arbitrary interval of time. We assume that assumptions {\bf (H1)}-{\bf (H2)} are satisfied as well as the following additional assumption on the autonomous equation for $X_m$.

{\bf (H3)} There exists an open interval $U\subset \bbr$ such that the SDE 
\begin{equation*}
dX_{m,t}=b_m(t,X_{m,t})dt \;+\sigma (X_{m,t})dW_t
\end{equation*}
possesses a unique strong solution taking values in $U,$ whenever $X_{m, 0} \in U.$ Moreover $\sigma(\cdot)$ is strictly positive on $U$ and its restriction to every compact interval in $U$ is of class $C^\infty$.

\subsection{Weak H\"ormander condition for (\ref{internal1})}\label{lwhforinternal1}
For convenience we rewrite system (\ref{internal1}) here:
\begin{eqnarray*}
dX_{1,t}&=& F (X_{l,t}, \ 1\leq l\leq m-1)dt +dX_{m,t}\, ,\\
dX_{i,t}&=&[-a_i(X_{1,t})X_{i,t}+b_i(X_{1,t})]dt, \, \, \, i = 2, \cdot, \cdot, \cdot , m-1\, ,\nonumber\\
dX_{m,t} &=& b_m(t,X_{m,t})dt \;+\sigma (X_{m,t})dW_t.\, \nonumber 
\end{eqnarray*}
The linearity of $dX_i$ , $i\in\{2,\cdot,\cdot,\cdot m-1 \}$ w.r.t $X_i$ has an important consequence that we recall in the following proposition (the proof of this proposition is provided in the Appendix). 
\begin{prop}\label{linear}
Fix $i\in\{2,\cdot\cdot\cdot,m-1\}$. Suppose that $X_{i,0}\in[0,1]$ a.s., and also $0\leq b_i(x)\leq a_i(x),$ for all $x\in \bbr$, $x$ denoting the first component of \eqref{internal1}. Then $\forall t>0,\, X_{i,t}\in[0,1]$ a.s.
\end{prop} 
In view of Proposition \ref{linear} we assume for the rest of this section that $X_{i,0}\in[0,1]$ a.s.,  and $0\leq b_i(x)\leq a_i(x),$ for all $x\in \bbr$, for all $i\in\{2,\cdot\cdot\cdot,m-1\}.$ We define $E_m:= \bbr\times[0,1]^{m-2}\times U$ where $U$ is given by {\bf (H3)}. Actually this is not a restriction since in the application we have in mind each variable $X_{i,u}$ describes the probability that some particular event occurs (see for instance section \ref{sec:main}).

\begin{defin}\label{determinant}
For any integer $k\geq 1$ denote by $\partial^{(k)} _{x_1}$ the partial derivative of order $k$ w.r.t. $x_1$. For any $x\in \bbr^{m-1} \times U \, $ consider $J_1(x):=F(x_1,x_2, \cdot,\cdot,\cdot,x_{m-1})$ and $J_i(x):=-a_i(x_1)x_i+b_i(x_1)$, $\ 2\leq i\leq m-1$. We define ${\bf D}(x)$ as the determinant of the matrix $(\partial^{(k)} _{x_1}J_i(x)\,; (i,k)\in \{ 1,\cdot, \cdot, \cdot, m-1\}^2)$.
\end{defin}
\begin{theo}\label{theo:internal_hoer}
Suppose that \eqref{internal1} satisfies {\bf (H1)}--{\bf (H3)}. Then
{\bf (LWH)} holds at any point $x=(x_i, \ 1\leq i\leq m)\in {\rm int} ( E_m) $ where ${\bf D}(x) \ne 0$. 
\end{theo}
The proof of Theorem \ref{theo:internal_hoer} will be given below (through Proposition \ref{crochetsdeliepour2}). First we state some important consequences and make some remarks. 
It is important to note that ${\bf D}(x)$ actually depends only on the $m-1$ first components of $x$. In particular if the $m-1$ first components of two points $x$ and $x'$ coincide, then ${\bf D}(x)={\bf D}(x')$. This remark will be important in the sequel (see e.g.\ Proposition \ref{cor:3} below). Moreover the condition in Theorem \ref{theo:internal_hoer} implies a version of {\bf (LWH)} uniform w.r.t. time on every compact interval $[0,T]$. Let us now define the set 
\begin{equation*}
{\cal D} := \{ (x_i, \ 1\leq i\leq m)\in {\rm int}(E_m); \quad{\bf D}(x) \ne 0\}  .
\end{equation*}
The set ${\cal D}$ is an open subset of $E_m$ by continuity of ${\bf D}$ on $\bbr^{m-1}\times U$. The following statement is a direct consequence of Theorems \ref{theo:main2} and \ref{theo:lsc} of Section \ref{localdensities}, taking $ E_m$ as state space.
\begin{theo}\label{localdensityinternal}
Suppose that \eqref{internal1} satisfies {\bf (H1)}--{\bf (H3)}.  Assume that $y_0\in {\cal D}$ and take $R>0$ such that $B_{3R}(y_0)\subset {\cal D}$. Then for any $x\in E_m$ and $t>0$, the random variable $X_t $ admits a density $ p_{ 0, t} ( x, \cdot ) $ on $ B_R (y_0 ) $. The map $ y \in B_R (y_0)\mapsto p_{ 0, t} ( x, y )$ is continuous, and for any fixed $ y \in B_R (y_0),$ the map $x \in E_m \mapsto  p_{0, t} (x,y)$ is lower semi-continuous.
\end{theo}

\begin{cor}\label{cor:111}
Grant the assumptions of Theorem \ref{localdensityinternal}. For all $x \in E_m , $ the following holds true. If there exists $y \in {\cal D} $ and $ t > 0 $ such that $ P_{0, t } ( x, U) > 0 $ for all sufficiently small neighborhoods $U$ of $y,$ then there exists $ \delta > 0 $ such that, if $K_1$ (resp.\ $ K_2$) denotes the closure of $ B_\delta ( x) $ (resp.\ $B_\delta ( y) $ ), 
$$ \inf_{ x' \in K_1 } \inf_{y' \in K_2} p_{0, t } (x', y') > 0 .$$
\end{cor}

The difficulty in practice is to obtain more information on ${\cal D}$,  in particular to know whether it coincides with ${\rm int}(E_m)$. At least one would like to be able to specify open regions included in ${\cal D}$. In general one can hope to achieve this goal only numerically unless the coefficients of the system are very simple. In Section \ref{sec:main} we provide details for a stochastic Hodgkin-Huxley model. The definition of ${\cal D}$ comes from a particular choice of successive Lie brackets where we look for the directions in space which propagate the noise at maximal possible speed according to the following intuition: the noise in (\ref{internal1}) is most rapidly transported through $X_1$ and $X_m$, since they are the only components carrying Brownian noise explicitly. Accordingly, except for the first Lie bracket $[A_0,A_1]$ which involves the drift $A_0$, we always use the diffusion coefficient $A_1$ in order to compute the brackets of higher order. The corresponding development of the solution of (\ref{internal1}) into iterated Ito integrals for small time steps $\delta$, shows that the speed of the diffusion is of order $ \delta^{\frac12}$ in the direction of $A_1$, of order $ \delta^{ 1 + \frac12}$ in the direction of $[A_0,A_1]$ and for the subsequent Lie brackets we add a factor $\frac{1}{2}$ to the exponent each time we use $A_1$, so that the speed of the diffusion is of order $\delta^2$ in the direction of $[[A_0,A_1],A_1]$, of order $\delta^{1 + 3 \times \frac12}$ in the direction of $[[[A_0,A_1],A_1],A_1]$ and so on. We refer the reader to \cite{NSW}, in particular identity (12).
Hence it is important to remember that belonging to ${\cal D}$ is only a sufficient condition for {\bf (LWH)} to hold. It may hold also at points outside ${\cal D}$ in which case the system suffers a slow down in the sense just explained.

\medskip

We now prove Theorem \ref{theo:internal_hoer} starting with the following key proposition about the computation of Lie brackets in this specific case. The proof is a direct consequence of the definition of Lie bracket recalled in section \ref{sec:hoerm} and is left to the reader. 
\begin{prop}\label{crochetsdeliepour2}
Consider on one hand $\varphi,\, \psi$ and $\rho$ smooth functions of $x_m$ defined on $U$ and on the other hand a family of smooth functions $y_i,\, 1\leq i\leq m-1,\, $ defined on $\bbr^{m-1},$ which do not depend on $x_m$. Let $\Xi$ and $Y$ denote vector fiels on $[0,+\infty]\times{\bbr}^m $of the following form, 
\begin{eqnarray*}
\Xi(t,x)&:=&\varphi(x_m)\, (\frac{\partial}{\partial x_1}+\frac{\partial}{\partial x_m}),\\
Y(t,x)&:=&\rho(x_m)\sum_{i=1}^{m-1}y_i\, \frac{\partial}{\partial x_i}+\psi(x_m)\,( \frac{\partial}{\partial x_1}+ \frac{\partial}{\partial x_m}).
\end{eqnarray*}
The Lie bracket $[\Xi,Y]$ takes the form
\begin{eqnarray*}
[\Xi,Y](t,x)&=&\varphi(x_m)\rho(x_m)\sum_{i=1}^{m-1}\, \partial_{x_1}y_i\, \, \frac{\partial}{\partial x_i}\\
&+&\varphi(x_m)\rho'(x_m)\sum_{i=1}^{m-1}\, y_i \frac{\partial}{\partial x_i}+(\varphi \psi'-\varphi' \psi)(x_m)\, \, (\frac{\partial}{\partial x_1}+\frac{\partial}{\partial x_m}).
\end{eqnarray*}
\end{prop}

\noindent {\bf Proof of Theorem \ref{theo:internal_hoer}.} According to the notations of Section \ref{sec:hoerm}, we write $A_1=\sigma(x_m)(\frac{\partial}{\partial x_1}+\frac{\partial}{\partial x_m})$ and $A_0=\frac{\partial}{\partial t}+\sum_{i=1}^m {\tilde b}_i\frac{\partial}{\partial x_i}$ where $\tilde b$ is given in (\ref{driftstratonovich}). Let us consider the Lie brackets defined recursively by $L_1:= [A_1,A_0]$ and $L_{k+1}=[A_1,L_k]$. In order to illustrate the relationship between the $L_k$ and the determinant ${\bf D}(x)$ introduced in Definition \ref{determinant}, we compute explicitly $L_1$ and $L_2$. We find first that
 \begin{equation*}
L_1=\sum_{i=1}^m\, \sigma(x_m)(\frac{\partial {\tilde b}_i}{\partial x_1}+\frac{\partial {\tilde b}_i}{\partial x_m})\frac{\partial}{\partial x_i}-\sigma'(x_m){\tilde b}_m(\frac{\partial}{\partial x_1}+\frac{\partial}{\partial x_m}). 
\end{equation*}
The drift $\tilde b$ in (\ref{internal1}) satisfies $\frac{\partial {\tilde b}_m}{\partial x_1}\equiv 0$, $\frac{\partial {\tilde b}_i}{\partial x_m}\equiv 0$ for all $i\in\{2,\cdot, \cdot, \cdot,m-1\}$. Moreover $\frac{\partial {\tilde b}_i}{\partial x_1}\equiv\partial_{x_1} J_i$ for all $i\in\{1,\cdot, \cdot, \cdot,m-1\}$. Hence
\begin{equation*}
L_1=\sum_{i=1}^{m-1}\, \sigma(x_m)\partial_{x_1} J_i\, \frac{\partial}{\partial x_i}+\sigma(x_m)(\frac{\partial {\tilde b}_1}{\partial x_m}\frac{\partial}{\partial x_1}+\frac{\partial {\tilde b}_m}{\partial x_m}\frac{\partial}{\partial x_m})-\sigma'(x_m){\tilde b}_m(\frac{\partial}{\partial x_1}+\frac{\partial}{\partial x_m}). 
\end{equation*}
We can further reduce the expression of $L_1$ using that the drift $\tilde b$ in (\ref{internal1}) satisfies also $\frac{\partial {\tilde b}_1}{\partial x_m}\equiv\frac{\partial {\tilde b}_m}{\partial x_m}$. We obtain
\begin{equation}
L_1=\sum_{i=1}^{m-1}\sigma(x_m)\partial_{x_1} J_i\, \frac{\partial}{\partial x_i}+\Big (\sigma(x_m)\frac{\partial {\tilde b}_m}{\partial x_m}-\sigma'(x_m){\tilde b}_m\Big )(\frac{\partial}{\partial x_1}+\frac{\partial}{\partial x_m}).
\end{equation}

Proposition \ref{crochetsdeliepour2} applies to $\Xi=A_1$ and $Y=L_1$ with $\varphi(x_m)\equiv\rho(x_m)\equiv \sigma(x_m),$ $y_i\equiv\partial_{x_1} J_i$, for $i\in\{1,\cdot, \cdot, \cdot,m-1\}$, $\psi(x_m)\equiv \sigma(x_m)\frac{\partial {\tilde b}_m}{\partial x_m}-\sigma'(x_m){\tilde b}_m$. Therefore, with this specific choice,
\begin{eqnarray}\label{secondbracket}
L_2&=&\sum_{i=1}^{m-1} \sigma(x_m)^2\, \partial_{x_1}^{(2)} J_i\,\,  \frac{\partial}{\partial x_i}\\
&+& \sum_{i=1}^{m-1} \sigma(x_m)\sigma'(x_m)\partial_{x_1} J_i\, \, \frac{\partial}{\partial x_i}+(\varphi \psi'-\varphi' \psi)(x_m)\, \, (\frac{\partial}{\partial x_1}+\frac{\partial}{\partial x_m})\nonumber.
\end{eqnarray}
Once again, identity (\ref{secondbracket}) coupled with Proposition \ref{crochetsdeliepour2} enables us to work by iteration. We thus obtain the following expression for $L_k$, for any $k\geq 1$:
\begin{equation}\label{kbracket}
L_k=\sum_{i=1}^{m-1}\sigma(x_m)^k\, \partial_{x_1}^{(k)} J_i\,\,  \frac{\partial}{\partial x_i}+\sum_{\ell=1}^{k-1}\sum_{i=1}^{m-1}\Phi_\ell(x_m)\, \partial_{x_1}^{(\ell)} J_i\, \, \frac{\partial}{\partial x_i}+\Phi(x_m)(\frac{\partial}{\partial x_1}+\frac{\partial}{\partial x_m}).
\end{equation}
The explicit expression of $\Phi_\ell$, $\Phi$ are not necessary to conclude. Indeed let us identify these Lie brackets with the column vectors in ${\bbr}^m$ obtained with their coordinates on the basis $(\frac{\partial}{\partial x_i}, \, i\in\{1,\cdot, \cdot, \cdot,m\})$. A sufficient condition for {\bf (LWH)} to be satisfied is that the vector space generated by $(A_1, L_k,\, k\in\{1,\cdot, \cdot, \cdot,m-1\}) $ coincides with ${\bbr}^m$. It is sufficient that the determinant formed with these vectors does not vanish. The definition of $A_1$ and formula (\ref{kbracket}) imply that this determinant coincides with the determinant obtained with the vectors $A_1$ and ${\tilde L}_k,\, k\in\{1,\cdot, \cdot, \cdot,m-1\}$ where ${\tilde L}_k:=\sum_{i=1}^{m-1}\sigma(x_m)^k\, \partial_{x_1}^{(k)} J_i\,\,  \frac{\partial}{\partial x_i}$. Since $\sigma$ does not vanish on $U$ (cf. {\bf (H3)}, we conclude that a sufficient condition is that ${\bf D}(x)$ does not vanish. \halmos

\subsection{Positivity of densities for models \eqref{internal1}}\label{positivedensities}
Once we have proved that densities exist for (\ref{internal1}), even if only locally, we look for regions where they are positive. For this purpose we combine control arguments and the support theorem. We keep the notation $E_m= \bbr\times[0,1]^{m-2}\times U$ introduced in the previous section. We start by proving an accessibility result for (\ref{internal1}) in Proposition \ref{theo:accessible} below, which holds without any assumption on the existence of densities and relies on some stability properties of the underlying deterministic system (\ref{subinternal}). We refer the reader to \cite{BLBMZ-13} for similar ideas in the framework of Piecewise Deterministic Markov Processes. 

Let $(X_u)_{u\geq 0}$ be a solution of (\ref{internal1}). We denote by $P_{0, t } ( x, \cdot )$ the law of $X_t$ when $X_0=x\,\,  {\rm a.s.}$ We recall (\ref{subinternal}) here for the reader's convenience:
\begin{eqnarray*}
dz_{1,t}&=& F (z_{l,t}, \ 1\leq l\leq m-1) dt \, , \\
dz_{i,t}&=&[-a_i(z_{1,t})z_{i,t}+b_i(z_{1,t})]dt, \, \, \, i = 2, \cdot, \cdot, \cdot , m-1.
 \end{eqnarray*}
\begin{prop}\label{theo:accessible}
Grant {\bf (H1)}--{\bf (H3)} and assume that $ U= \bbr .$ 
We keep moreover the assumptions of Proposition \ref{linear} and suppose that $0 < b_i < a_i$ for all $i \in\{2, \cdot, \cdot, \cdot , m-1\}$. Given an arbitrary real number $z_1$, consider $\overline z:= (z_1,\,  y_{i,\infty}(z_1),\, \,i \in\{2, \cdot, \cdot, \cdot , m-1\})$ in $\bbr^{m-1}$, where $y_{i,\infty}(z_1):=\frac{b_i(z_1)}{a_i(z_1)}$ is an equilibrium point for the $i-$th equation when we keep the first variable fixed at constant value $z_1.$ For all $x\in E_m$ and any neighborhood ${\cal N}$ of $ \overline z$ in $\bbr\times(0,1)^{m-2}$ there exists $t_0 $ such that  
\begin{equation}
 \forall t \geq t_0,\quad P_{0, t } ( x ,\, {\cal N} \times \bbr ) > 0 .
\end{equation}
\end{prop}

\begin{prop}\label{cor:3.1} 
Let us keep the notations and assumptions of Proposition \ref{theo:accessible}. Consider an arbitrary real number $z_1$ and the associated point $\overline z:= (z_1,\,  y_{i,\infty}(z_1),\, \,i \in\{2, \cdot, \cdot, \cdot , m-1\})$ in $\bbr^{m-1}$. Assume that ${\bf D} (\overline z,u)\neq 0$ for some $u\in \bbr$, where the determinant ${\bf D}$ has been introduced in Definition \ref{determinant}. Then for all $x\in E_m$, there is $t_0 > 0 $ such that for all $t \geq t_0 $ the following holds true. There exist $ \overline u = \overline u (t) \in \bbr $ and  $\delta = \delta ( t) >0$ such that, if $ K_1$ (resp. $K_2$) denotes the closure of $B_\delta (x)$ (resp. $B_\delta (\overline z,\overline u)$),
\begin{equation*}
 \inf_{ x' \in K_1 } \inf_{y'\in K_2}  p_{ 0, t} (x', y' )> 0.
\end{equation*}
\end{prop}
Notice that for each $i \in\{2, \cdot, \cdot, \cdot , m-1\}$, the solution of $dy_t = (-a_i(z_1) y_t+b_i(z_1))dt$ with $z_1$ as a fixed parameter, converges to $y_{i,\infty} (z_1)$ when $t\rightarrow +\infty$ and that $y_{i,\infty} (z_1)$ is globally asymptotically stable. Proposition \ref{theo:accessible} holds in particular when $F(\overline z)=0$. In this case $\overline z$ is an equilibrium point of (\ref{subinternal}).

\noindent{\bf Proof of Proposition \ref{theo:accessible}.}\label{sec:supp} Let $z_1\in\bbr$ and the associated point $\overline z:= (z_1,\,  y_{i,\infty}(z_1),\, \,i \in\{2, \cdot, \cdot, \cdot , m-1\})$ in $\bbr^{m-1}$. As in the proof of Theorem \ref{theo:lsc} we write 
$ \Omega$ for $C ( [ 0, \infty [ , \bbr^m )  $
and endow it with its canonical filtration $ ( {\cal F}_t)_{t \geq 0 } .$ Recall that $ \mathbb{P}_{0,x}$ is the law of $(X_u)_{u \geq 0}$ starting from $x$ at time $0.$ We first localize the system by a sequence of compacts $(K_n)$ according to {\bf (H1)} and let $T_n = \inf \{ t : X_t \in K_n^c \} $ be the exit time of $K_n$. For a fixed $n,$ let $ b^{(n)} (t, x) $ and $\sigma^{(n)} (x)$ be $C_b^\infty -$extensions in $x$ of $ b(t, \cdot_{| K_n }) $ and $\sigma_{| K_n}$ respectively and $X^{(n)}$ be the associated diffusion process (here we denote the coefficients of (\ref{internal1}) by $b$ and $\sigma$ for short). For any integer $n\geq 1$ and starting point $x,$ we write $\mathbb{P}_{0, x}^{(n )}$ for the law of $(X^{(n)}_u)_{u \geq0}$ on $\Omega$ satisfying $X_0^{(n)}=x$. We wish to find lower bounds for quantities of the form $\mathbb{P}_{0,  x } (  B ) $ where $B = \{ f \in \Omega : f(t) \in {\cal N}\times \bbr \}  \in {\cal F}_t$, for any $t>0$ given . We start with the following inequality which holds for any $t > 0 $ and $n$:
\begin{equation}\label{eq:tobelb}
\mathbb{P}_{0,  x } (  B ) \geq \mathbb{P}_{ 0, x } ( \{ f  \in B ; T_n > t \}  ) 
 = \mathbb{P}^{(n)}_{ 0, x } (  \{ f \in B ; T_n >  t \} )  .
\end{equation}
In the sequel we show that for some integer $n_0$ and any fixed $x \in K_{n_0}$, the quantity $\mathbb{P}^{(n)}_{ 0, x } (  \{ f \in B ; T_n >  t \} )$ is indeed positive provided that $n$ is sufficiently large. We are therefore interested in the support of $\mathbb{P}^{(n)}_{ 0, x }$. Fix $t$ and let $ {\cal C} := \{ {\tt h} : [ 0, t ] \to \bbr : {\tt h}(s) = \int_0^s \dot {\tt h} (u) du,\,  \forall \, s \le t, \, \int_0^t \dot {\tt h}^2 (u) du < \infty \} $ be the Cameron-Martin space. Given ${\tt h} \in {\cal C} ,$ consider $ X({\tt h})\in \bbr^m$ the solution of the differential equation 
\begin{equation}\label{eqcontrol1}
X({\tt h})_s = x + \int_0^s \sigma^{(n)} ( X({\tt h})_u) \dot {\tt h} (u) du + \int_0^s  {\tilde b}^{(n)} ( u, X({\tt h})_u) du ,\quad s\le t.
\end{equation}
If (\ref{eqcontrol1}) were time homogeneous, the support theorem would imply that the support of $\mathbb{P}^{(n)}_{0, x}$ in restriction to $ {\cal F}_t$ is the closure of the set $ \{ X({\tt h}) : {\tt h} \in {\cal C} \} $ with respect to the uniform norm on $ [ 0, t]$ (see e.g.  \cite{Millet-Sanz} Theorem 3.5 or \cite{BenArous} Theorem 4). To conclude in our situation as well, it is enough to replace the $m-$dimensional process $X^{(n)}$ by the $(m+1)-$dimensional process $(t, X^{(n)} _t) $ which is time-homogenous. In order to proceed further we construct a control ${\tt h}$ so that $X({\tt h})$ remains in $K_n$ during $ [ 0, t ]$ provided $n$ is sufficiently large. We start by exploiting stability properties of the underlying deterministic system (\ref{subinternal}). Let $ \gamma: \bbr\mapsto\bbr$ a smooth function satisfying $\gamma(\tau):=z _1$ for all $\tau\geq 1$. Consider $y_s\in\bbr^{m-2}$ solving
\begin{equation*}
dy_{i,s}=[-a_i(\gamma(s))y_{i,s}+b_i(\gamma(s))]ds, \, \, \, i = 2, \cdot, \cdot, \cdot , m-1.
\end{equation*}
Then for all $t>1$, 
\begin{eqnarray*}
y _{i,t} &=& y_{i,0} e^{ - \int_0^t  a_i ( \gamma(s)) ds } + \int_0^t b_i ( \gamma(u)) e^{ - \int_u^t a_i( \gamma(r)) dr } du\\
&=&y_{i,1}e^{ - a_i ( z_1)(t-1)}+y_{i,\infty}(z_1)(1-e^{ - a_i ( z_1)(t-1)}).
\end{eqnarray*}
This formula expresses the fact that on $[1,+\infty[$, the coefficients $a_i(\gamma(s))$  (resp. $b_i(\gamma(s))$) are constant equal to $a_i ( z_1)$ (resp. $b_i ( z_1)$). Hence for any  $\varepsilon>0$ there exists $t_0>1$ such that $|y_{i,t}-y_{i,\infty}(z_1)|<\varepsilon$ for all $t\geq t_0$ and all $2\leq i\leq m-2$. Now take $\varepsilon$ so small that $ B_\varepsilon ( \overline z) \subset {\cal N} .$ Then for all $t\geq t_0>1, $ the vector $(\gamma(t), y _{i,t}\, 2\leq i\leq m-2)$ belongs to $ B_\varepsilon ( \overline z ) $ (remember that for $t>1$, $\gamma(t)$ is fixed at $z_1$). \\ 
Fix an integer $n_0$ and $x$ in $K_{n_0}$. We are now able to construct a control $h\in{\cal C}$ such that the solution of (\ref{eqcontrol1}) remains in $K_n$ during finite time intervals for all $n$ large enough. Choose a function $\gamma$ as above satisfying moreover $ \gamma (0 ) = x_1$, $ \gamma (1 ) =z _1$. Define $(Z_s)_{s\geq 0}\in \bbr^m$, the deterministic path starting from $x$ such that 
\begin{eqnarray}\label{construction}
Z_{1,s}&=&\gamma(s)\, ,  \\
dZ_{i,s}&=&[-a_i(Z_{1,s})Z_{i,s}+b_i(Z_{1,s})]ds, \, \, \, i = 2, \cdot, \cdot, \cdot , m-1\, ,\nonumber\\
Z_{m,s}&=& x_m-x_1+\gamma(s)-\int_0^s F(Z_u)du.\nonumber
 \end{eqnarray}
Note that  $(Z_s,\, s\in[0,t])$ is bounded and therefore remains in $K_n$ for all $n$ large enough. 
Now fix $t \geq t_0 $ and consider  a function ${\tt h}$ defined by
\begin{equation}\label{hdot}
 \dot {\tt h} (s) := \frac{\dot\gamma(s)- F(Z_s )-b_m(s,Z_{m,s})+\frac{1}{2}\sigma(Z_{m,s}) \sigma'(Z_{m,s})} {\sigma(Z_{m,s})}.
\end{equation}
Since by assumption $\sigma (\cdot ) > 0 $ on $ \bbr $ (we have assumed that $ U = \bbr $), the expression (\ref{hdot}) is well-defined. This assumption also provides that $\dot{\tt h} \in L^2 ( [ 0, t ] ) ,$ hence ${\tt h} \in {\cal C} .$ Hence, with such a choice of ${\tt h},$ the solution $X({\tt h})$ of equation (\ref{eqcontrol1}) coincides with the solution $Z$ of system (\ref{construction}).  As explained previously, we can choose $n$ such that $(Z_s,\, s\in[0,t])$ remains in $K_n$. \\
Consider now, for $\delta>0$, the tubular neighborhood $T_\delta$ of $(Z_s,\, s\in[0,t])$ in $\Omega$ of size $\delta$, namely the set  $\{ f \in  \Omega  : \sup_{ s \le t } | f(s) - Z_s | <  \delta \} .$ By the support theorem $ \mathbb{P}^{(n)}_{0, x} (T_\delta) > 0 .$ Remember that we have chosen $\epsilon$ and $t_0$ in order to satisfy $T_\delta \subset   \{ f \in \Omega : f(t) \in B_\varepsilon (\overline z) \times \bbr \}    $ as well as $B_\varepsilon (\overline z)\subset{\cal N}$. Choosing $ \delta \le \varepsilon/2  $ such that $ T_\delta\subset \{ f \in \Omega : T_n (f) > t \} $, we conclude as announced that 
\begin{equation*}
 P_{0,t}(x,{\cal N}\times \bbr)\geq  P_x ( X_t \in B_\varepsilon (\overline z) \times \bbr  ) \geq \mathbb{P}^{(n)}_{0, x } (T_\delta) > 0 . 
 \end{equation*}
\halmos



\noindent{\bf Proof of Proposition \ref{cor:3.1}.}\label{sec:supp3} The fact that ${\bf D} (\overline z,u)\neq 0$ for some $ u \in \bbr $ implies that ${\bf D} (\overline z, \tilde u)\neq 0$  for all $ \tilde u \in \bbr ,$ by Theorem \ref{theo:internal_hoer}. The attainability at time $t$ is proven as in the proof of Proposition \ref{theo:accessible}. For $x \in E_m, $ $\gamma_1 $ and $t_0$ as there, $t_0 > 1 ,$ we define for $ t \geq t_0$ 
$$ \overline u(t) = x_m - x_1 + z_1 - \int_0^t F ( Z_s) ds \in U = \bbr .$$ 
Then there is some $ \delta (t) > 0 $ such that 
$$ P_{0, t } ( x, B_{\delta (t) } ( \overline{z}, \overline u(t) ) ) > 0.$$
Applying Corollary \ref{cor:111} to $y = (  \overline{z}, \overline u(t) ) \in {\cal D} $ 
finishes the proof. \halmos

We now show that, during any arbitrary long period, with positive probability, the stochastic system (\ref{internal1}) is able to reproduce the behavior of $(d z_t,I(t))\in \bbr^m$ where \begin{eqnarray}\label{subinternalwithI}
dz_{1,t}&=& [F (z_{l,t}, \ 1\leq l\leq m-1) +I(t)]dt  \\
dz_{i,t}&=&[-a_i(z_{1,t})z_{i,t}+b_i(z_{1,t})]dt, \, \, \, i = 2, \cdot, \cdot, \cdot , m-1 \nonumber
 \end{eqnarray}
with $I(t)$ an arbitrary smooth input applied to (\ref{subinternal}). Note that by comparing (\ref{internal1}) and (\ref{subinternalwithI}) we see that the $m-$the component $X_m$ of the stochastic system \eqref{internal1} has to be compared to a deterministic control path (\ref{subinternalwithI}) to which we add an $m-$th coordinate given by $ t \to X_{m, 0 } + \int_0^t I(s) ds .$

Remember $ \, B_{\delta}( x)$ denotes the open ball of radius $\delta$ centered at $x$. In the following two propositions, $U$ is again an open interval in $ \bbr.$ 
\begin{prop}\label{theo:posinternal}
Suppose that \eqref{internal1} satisfies {\bf (H1)}--{\bf (H3)}. 
Fix $x\in E_m$ and $t>0$. Let $I$ be a smooth deterministic input such that $x_m+\int_0^s I(r)dr\in U$ for all $s\leq t$. Define $\mathbb{X}_s^x:=(\mathbb{Y}_s^{\tilde x}, \, x_m+\int_0^s I(r)dr, \, s\leq t)$ where $\mathbb{Y}^{\tilde x}$ is the deterministic path solution of (\ref{subinternalwithI}) starting from $\tilde x:=(x_i, 1\leq i\leq m-1)$. We denote by $\mathbb{P}_{0, x}$ the law of the solution of (\ref{internal1}) starting at $x$. Then for any $\varepsilon > 0$ 
\begin{equation*}
\mathbb{P}_{ 0,x} \left( \left\{ f \in \Omega : \sup_{ s \le t } | f(s) - \mathbb{X}^x_s | \le \varepsilon \right\} \right) > 0 
\end{equation*}
and moreover there exists $\delta > 0$ such that for all $ x'' \in B_\delta ( x)$ 
\begin{equation*}
 \mathbb{P}_{ 0,x''} \left( \left\{ f \in \Omega : \sup_{ s \le t } | f(s) - \mathbb{X}^{x }_s | \le \varepsilon \right\} \right) > 0 .
 \end{equation*}
\end{prop}

\noindent{\bf Proof of Proposition \ref{theo:posinternal}.}\label{sec:supp2} We keep the notations introduced in the proof of Proposition \ref{theo:accessible}. In the course of this proof we have shown that the support theorem applies to inhomogeneous diffusions like the one obtained after localizing (\ref{internal1}). Moreover we still hope to reach the positivity we are looking for through inequalities (\ref{eq:tobelb}) and paths solving (\ref{eqcontrol1}) for $h\in {\cal C}$, that remain in $K_n$ during $ [0, t]$ for $n$ sufficiently large.  So the system we work with is the localized one. Consider $I$ to be a deterministic input such that $x_m+\int_0^s I(r)dr\in U$ for all $s\leq t$. Define $\chi_{m,s}:=x_m+\int_0^s I(r)dr$ for all $s\leq t$ and 
\begin{equation}\label{hdot2}
 \dot {\tt h} (s) := \frac{I(s)-b_m(s,\chi_{m,s})+\frac{1}{2}\sigma(\chi_{m,s})\sigma'(\chi_{m,s})}{\sigma(\chi_{m,s})}.
\end{equation}
By definition  $(\chi_{m,s}, s\leq t)$ lies in a compact interval included in $U$. Then, the expression (\ref{hdot2}) is well-defined by assumption {\bf(H3)}. This assumption also provides that $\dot{\tt h} \in L^2 ( [ 0, t] )$ hence ${\tt h} \in {\cal C}$. Moreover, with such a choice of ${\tt h},$ the controlled path $X({\tt h})$, solution of (\ref{eqcontrol1}), coincides with $(\mathbb{Y}_s^{\tilde x}, \, \chi_{m,s}, \, s\leq t)$ where $\mathbb{Y}^{\tilde x}$ is the deterministic path solution of (\ref{subinternalwithI}) starting from $\tilde x=(x_i, 1\leq i\leq m-1)$. We can choose $n$ large enough such that $(\mathbb{Y}_s^{\tilde x}, \, \chi_{m,s},\, s\in[0,t])$ remains in $K_n$. We write $\mathbb{X}_s^x$ for $(\mathbb{Y}_s^{\tilde x}, \, \chi_{m,s})$. 
Remember that $\Omega=C([0,\infty[;\bbr^m)$ and for $\delta>0$, consider the tubular neighborhood $T_\delta$ of $\mathbb{X}^x$ on $[0,t]$ namely the set  $\{ f \in  \Omega  : \sup_{ s \le t } | f(s) - \mathbb{X}^x_s | <  \delta \} .$ By the support theorem $ \mathbb{P}^{(n)}_{0, x} (T_\delta) > 0.$ Choose now $ \delta$ such that $ T_\delta\subset \{ f \in \Omega : T_n (f) > t \} $. Taking $T_\delta$ as the set $B$ in (\ref{eq:tobelb}) yields the first statement of Proposition \ref{theo:posinternal}. The second one follows from the Feller property of $\mathbb{P}^{(n)}_{0, x }$ which enables us to extend the first statement to a small ball around $x$.
\halmos

We close this section with the following consequence of Proposition \ref{theo:posinternal} from which we borrow the notations. In addition, we assume that the deterministic system \eqref{subinternal} admits equilibria as considered in \eqref{eq:5}.

\begin{prop}\label{cor:3} 
Assume that \eqref{internal1} satisfies {\bf (H1)}--{\bf (H3)}, that $z^*$ is an equilibrium point of (\ref{subinternal}) such that $ z^* \in \bbr \times (0, 1 )^3$ and that $ r \to I(r) $ is some smooth input. Consider $x_m \in U$ and $t$ such that $x_m+\int_0^s I(r)dr \in U $ for all $ 0 \le s \le t.$ Let  $x:=(z^*,x_m)$ and $y :=(z^*, \, x_m+\int_0^t I(r)dr).$ Then the following holds. If ${\bf D} (x)\neq 0,$ then there exists $\delta>0$ such that  
$$ \inf_{ x'\in K } \inf_{y' \in K'}  p_{ 0, t} (x', y') > 0 ,$$
where $ K$ (resp. $K'$) stands for the closure of $B_\delta (x)$ (resp. $B_\delta (y)$). 
\end{prop}

\noindent{\bf Proof of Proposition \ref{cor:3}.} From $ {\bf D } (x) \neq 0 $ we have $  {\bf D } (y) \neq 0 ,$ by Theorem \ref{theo:internal_hoer}. Moreover, $P_{0, t } ( x, B_\delta (y)) > 0 ,$ as in the proof of Proposition \ref{theo:posinternal}. Applying Corollary \ref{cor:111} to $ y \in {\cal D} $ finishes the proof. 
\halmos

\section{Application to physiology}\label{sec:main}
In this section we apply the above results to a random system based on the Hodgkin-Huxley model well known in physiology. This random system belongs to the family of SDEs with internal variables and random input presented in section \ref{ivri}. We start by some reminders on the deterministic Hodgkin-Huxley model that we call (HH) for short. 

\subsection{The deterministic (HH) system}\label{sec:1.0} 

The deterministic Hodgkin-Huxley model for the membrane potential of a neuron (cf \cite{HH-52}) has been extensively studied over the last decades. There seems to be a large agreement that it models adequately many observations made on the response to an external input, in many types of neurons. This model belongs to the family of conductance-based models. Indeed it features two types of voltage-gated ion channels responsible for the import of Na$^+$ and export of K$^+$ ions through the membrane. The time dependent conductance of a sodium (resp. potassium) channel depends on the state of four gates which can be open or closed; it is maximal when all gates are open. There are two types of gates $m$ and $h$ for sodium, one type $n$ for potassium. The variables $n_t$, $m_t$, $h_t$ describe the probability that a gate of corresponding type be open at time $t$. 
Then, the Hodgkin-Huxley equations with deterministic input $I$ which may be time dependent, is the $4D$ system 
\begin{eqnarray}\label{HH}
dV_t  \;&=&\; I(t)\, dt \;- \left[\, \ov g_{\rm K}\,n_t^4\, (V_t-E_{\rm K}) \;+\; \ov g_{\rm Na}\,m_t^3\, h_t\, (V_t  -E_{\rm Na}) \;+\; \ov g_{\rm L}\, (V_t-E_{\rm L}) \right] dt,\\
dn_t \;&=&\;  \left[\, \al_n(V_t)\,(1-n_t)  \;-\; \beta_n(V_t)\, n_t  \,\right] dt, \nonumber \\
dm_t \;&=&\;  \left[\, \al_m(V_t)\,(1-m_t)  \;-\; \beta_m(V_t)\, m_t  \,\right] dt, \nonumber  \\
dh_t \;&=&\;  \left[\, \al_h(V_t)\,(1-h_t)  \;-\; \beta_h(V_t)\, h_t  \,\right] dt ,\nonumber 
\end{eqnarray}
where we adopt the notations and constants of \cite{I-09}. For instance the conductance of a sodium channel at time $t$ is given by $\ov g_{\rm Na}\,m_t^3\, h_t$. The functions $\al_n, \beta_n, \al_m, \beta_m, \al_h, \beta_h$ take values in $(0,\infty)$ and are analytic, i.e. they admit a power series representation on $\bbr$. They are given as follows:
\begin{equation}
\begin{array}{llllll}
\alpha_n(v)  &=& \frac{0.1-0.01v }{\exp(1-0.1v)-1}, &  \beta_n(v) &= &0.125\exp(-v/80) ,  \\
\alpha_m(v)& = &\frac{2.5-0.1v}{\exp(2.5-0.1v)-1} , & \beta_m(v)&= &4\exp(-v/18) ,  \\
\alpha_h(v) &= &0.07\exp(-v/20) , &\beta_h(v) &=& \frac{1}{\exp(3-0.1v)+1}.
\end{array}
\end{equation}
Moreover if we set $a_n:=\al_n+\beta_n$, $b_n:=\al_n$ and analogously for $m$ and $h$, we see that (HH) can be written as a particular case of (\ref{subinternal}) with $F$ given by
\begin{equation}\label{eq:F}
F(v,n,m,h) =-[ \ov g_{\rm K}\,n^4\, (v-E_{\rm K})+\ov g_{\rm Na}\,m^3\, h\, (v  -E_{\rm Na}) \;+\ov g_{\rm L}\, (v-E_{\rm L})]. 
\end{equation}
If the variable $V$ is kept constant at $v\in\bbr$, the variables $n_t$, $m_t$, $h_t$ converge when $t\rightarrow +\infty$ respectively towards
\begin{equation}\label{eq:ninfty}
n_\infty(v) := \frac{\al_n}{\al_n+\beta_n}(v) \;,\; m_\infty(v) := \frac{\al_m}{\al_m+\beta_m}(v) \;,\; h_\infty(v) := \frac{\al_h}{\al_h+\beta_h}(v) \;.
\end{equation}
The parameter $\ov g_{\rm Na}$ (resp. $\ov g_{\rm K}$) is the maximal conductance of a sodium (resp. potassium) channel while $\ov g_{\rm L} $ is the leak conductance. The parameters $E_{\rm K} $, $E_{\rm Na} $, $E_{\rm L}$ are called reversal potentials. Their values $\ov g_{\rm K} = 36$, $\ov g_{\rm Na} = 120$, $\ov g_{\rm L} = 0. 3$ , $E_{\rm K} = - 12$, $E_{\rm Na} = 120$, $E_{\rm L} =  10.6$ are those of \cite{I-09}.  

The Hodgkin-Huxley system exhibits a broad range of possible and qualitatively quite different behaviors, depending on the specific input $I$. In response to a periodic input, the solution of (\ref{HH}) displays a periodic behavior (regular spiking of the neuron on a long time window) only in special situations. Let us first mention that there exists some interval $U$ such that time-constant input in $U$ results in periodic behavior for the solution of (\ref{HH}) (see \cite{RM-80}). For an oscillating input, there exists some interval $J$ such that oscillating inputs with frequencies in $J$ yield periodic behavior (see \cite{AMI-84}). Periodic behavior includes that the period of the output can be a multiple of the period of the input. However, the input frequency has to be compatible with a range of preferred frequencies of (\ref{HH}), a fact which is similarly encountered in biological observations (see \cite{I-09}). Indeed there are also intervals $\wt I$ and $\wt J$ for which time-constant input in $\wt I$ or oscillating input at frequency $f\in \wt J$ leads to chaotic behavior. Using numerical methods \cite{Endler} gives a complete tableau. 

\subsection{(\ref{HH}) with random input}
It has been shown in \cite{PTW-10} that conductance-based models like (\ref{HH}) are fluid limits of a sequence of Piecewise Deterministic Markov Processes. Such limit theorems enable to study the impact of {\it channel noise} (also  called {\it intrinsic noise}) on latency coding. Our setting is different. The noise here is external coming from the network in which the neuron is embedded, through its dendritic system. This system has a complicated topological structure and carries a large number of synapses which register spike trains emitted from a large number of other neurons within the same active network. We model the cumulated dendritic input as a diffusion of mean-reverting type carrying a deterministic signal $S$. The resulting system that we consider is the following particular case of (\ref{internal1}):
\begin{eqnarray}\label{HHrandominput}
dV_t  \;&=&\; d\xi_t \;- \left[\, \ov g_{\rm K}\,n_t^4\, (V_t-E_{\rm K}) \;+\; \ov g_{\rm Na}\,m_t^3\, h_t\, (V_t  -E_{\rm Na}) \;+\; \ov g_{\rm L}\, (V_t-E_{\rm L}) \right] dt\, ,\\
dn_t \;&=&\;  \left[\, \al_n(V_t)\,(1-n_t)  \;-\; \beta_n(V_t)\, n_t  \,\right] dt \, ,\nonumber \\
dm_t \;&=&\;  \left[\, \al_m(V_t)\,(1-m_t)  \;-\; \beta_m(V_t)\, m_t  \,\right] dt \, , \nonumber\\
dh_t \;&=&\;  \left[\, \al_h(V_t)\,(1-h_t)  \;-\; \beta_h(V_t)\, h_t  \,\right] dt\, ,\nonumber   \\
d\xi_t \;&=&\; (\, S(t)-\xi_t\,)\, \tau dt \;+\; \gamma\, q (\xi_t)\, \sqrt{\tau} dW_t \, ,\nonumber
\end{eqnarray}
parametrized in terms of $\tau$ (governing speed) and $\gamma$ (governing spread). For instance $\xi$ can be of Ornstein-Uhlenbeck (OU) type (then $U=\bbr$, $q (\cdot)\equiv 1$) or of Cox-Ingersoll-Ross (CIR) type (then $U=(-K,\infty)$, $q (x)=\sqrt{(x+K)\vee 0\;}$ for $x\in U$, and $K$ is chosen in $]\frac{\gamma^2}{2}+\sup|S|,+\infty[$). Such a choice builds on the statistical study \cite{H-07}. When the deterministic signal $S$ is periodic, it is shown in \cite{HK-10} that $\xi $ of  OU type admits a periodically invariant regime under which the signal $S(\cdot) $ is related to expectations of $\xi$ via the formula $s \to E_{\pi, 0} ( \xi_s ) =  \int_0^\infty S(s - \frac{r}{\tau} ) e^{ - r } dr $. In the companion papers \cite{HLT-2} and \cite{HLT-3}  we address the periodic ergodicity of the solution to (\ref{HHrandominput}). Ergodicity properties when $\xi$ is of OU type are the topic of \cite{HLT-2}. The case of CIR is covered in \cite{HLT-3} where also limit theorems are proved. Below we will conduct a numerical study of {\bf (LWH}) for (\ref{HHrandominput}), based on Theorem \ref{theo:internal_hoer}. In this theorem, the specific nature of $\xi$ plays no role in the definition of the determinant ${\bf D}$ provided that the SDE satisfied by $\xi$ satisfies assumption {\bf (H3)}, cf.\ Proposition \ref{weakhoerHH} below. Therefore, the results of this numerical study apply to general random (HH) where we replace the last line in \eqref{HHrandominput} by $d \xi_t = b_5 (t , \xi_t) dt + \sigma ( \xi_t) d W_t$.

\subsection{Weak H\"ormander condition for (\ref{HHrandominput})}\label{secweakhoerHH}
\subsubsection{The determinant ${\bf \Delta}$}
Applying Theorem \ref{theo:internal_hoer} and Definition \ref{determinant} we have to consider points where the $4D$ determinant, whose columns are the partial derivatives of the coefficients of (\ref{HH}) with respect to the first variable $v$ from order one to order four, does not vanish. Since in this case the function $F$ given in (\ref{eq:F}) is linear in $v$, we obtain that $\partial_v^{(k)}F=0$ for $k\in\{2,3,4\}$. Moreover $\partial_v F(v,n,m,h)=-( \ov g_{\rm K}\,n^4+\ov g_{\rm Na}\,m^3\, h+\ov g_{\rm L})$ never vanishes on $[0,1]^3$. So actually in this case, it is sufficient to consider a $3D$ determinant extracted from ${\bf D}$.\begin{prop}\label{weakhoerHH}
Assume that $\sigma$ remains strictly positive on $U$. Let us introduce the notation $d_n(v,n) :=  -a_n(v)n+b_n(v)$ and analogous ones for $m$ and $h$. Then {\bf (LWH)} for (\ref{HHrandominput}) is satisfied at any point $(v,n,m,h,\zeta)\in \bbr\times(0,1)^3\times U$ where ${\Delta}(v,n,m,h)\neq 0$ with 
\begin{equation}\label{eq:det}
{\bf \Delta}(v,n,m,h) \;:=\; \det \left( \begin{array}{lll} 
\partial_v^{(2)} d_n  & \partial_v^{(3)} d_n & \partial_v^{(4)} d_n \\
\partial_v^{(2)} d_m & \partial_v^{(3)} d_m &  \partial_v^{(4)} d_m  \\
\partial_v^{(2)} d_h & \partial_v^{(3)} d_h & \partial_v^{(4)} d_h  \\
\end{array} \right). 
\end{equation}
\end{prop}

\begin{prop}
The set of points in $(v,n,m,h,\zeta)\in \bbr\times(0,1)^3\times U$ where ${\bf \Delta}$ does not vanish has full Lebesgue measure.
\end{prop}

{\bf Proof.}
We say that a set has full Lebesgue measure if its complement has Lebesgue measure zero. Firstly it can be shown numerically that indeed there exists points $(v,n,m,h,\zeta)$ such that  ${\bf \Delta}(v,n,m,h)\neq 0$ (see Section \ref{ex:3} below). Moreover, for any fixed $v \in \bbr ,$ the function $ (n,m,h) \mapsto {\bf \Delta}( v, n,m,h) $ is a polynomial of degree three in the variables $n,m,h .$ In particular, for any fixed $v,$ either ${\bf \Delta}(v,.,.,.)$ vanishes identically on $(0,1)^3$, or its zeros form a two-dimensional sub-manifold of $(0,1)^3 .$ 
Finally, since ${\bf \Delta}$ is a sum of terms
$$
\mbox{(some power series in $v$)} \cdot n^{\vep_n} m^{\vep_m} h^{\vep_h}
$$
with epsilons taking values 0 or 1, it is impossible to have small open $v$-intervals where it vanishes identically on $(0,1)^3 .$ We conclude the proof by integrating over $v$ and using Fubini's Theorem. \halmos

\noindent Although the condition ${\bf \Delta}\neq 0$ is only a sufficient condition ensuring that {\bf (LWH)} is satisfied locally, it is convenient since it is possible to evaluate ${\bf \Delta}(v,n,m,h)$ numerically. This is done in section \ref{ex:3} below. However we are not able to characterize the whole set of points where {\bf (LWH)} holds unless we make more stringent assumptions on the input $\xi$, like for instance to assume that the coefficients of the SDE it satisfies are analytic (cf. \cite{HLT-2}, \cite{HLT-3}).

\subsection{Numerical study of the determinant ${\bf \Delta}$}\label{ex:3}
We compute numerically the value of ${\bf \Delta}$ at points of the form $(v, n_\infty(v),m_\infty(v), h_\infty(v))$ as in (\ref{eq:ninfty}). The function $F_\infty(v):=F(v, n_\infty(v),m_\infty(v), h_\infty(v))$ is strictly increasing on an interval ${\cal I}$ containing ${\cal I}_0=(- 15, + 30)$ hence it defines a bijection between the constant input $I(t)=c$ in (\ref{HH}) and the solution of the equation $F_\infty(v)=c$ that we denote by $v_c$. Therefore for any $v\in {\cal I}$, the point $(v, n_\infty(v),m_\infty(v), h_\infty(v))$ is the equilibrium point of (\ref{HH}) submitted to the constant input $c=F_\infty(v)$. We use this fact below since it may be more convenient to work with $v$ than with $c$ even if classically one considers $c$ as the parameter of interest. For instance the point $(0, n_\infty(0), m_\infty(0), h_\infty(0) )$ corresponds to $c=F( 0, n_\infty(0), m_\infty(0), h_\infty(0) ) \approx -0.0534$. We found that ${\bf\Delta}(0, n_\infty(0), m_\infty(0), h_\infty(0) )<0$ and moreover the function $v \mapsto {\bf \Delta}\left(v, n_\infty(v), m_\infty(v), h_\infty(v)\right)$ has exactly two zeros on the interval ${\cal I}_0=(- 15, + 30)$ located at $v \approx -11.4796$ and $v \approx +10.3444$. Hence for all values of $c$ belonging to $]F_\infty ( -10) , F_\infty  (+10)[= ]-6.15, 26.61[ $, the determinant ${\bf \Delta}(v_c, n_\infty(v_c),m_\infty(v_c), h_\infty(v_c))$ remains strictly negative.
   
\noindent For constant input $c=15$, the equilibrium point $(v_c, n_\infty(v_c),m_\infty(v_c), h_\infty(v_c))$ is unstable and (\ref{HH}) possesses a stable orbit (see Figure 1 where we plotted $n$ against $v$ the long of the orbit). We studied $ t \mapsto {\bf \Delta}(v_t,n_t,m_t,h_t)$ along this orbit. The periodic behavior is displayed in Figure 3, starting with a numerical approximation of the unstable equilibrium point and showing that the system switches towards a stable orbit. In this picture, already the last four orbits can be superposed almost perfectly. The value of ${\bf\Delta}$ at equidistant time epochs on the last complete orbit (starting and ending when the membrane potential $v$ up-crosses the level $0$, and having its spike near time $t=180$) is provided in Figure 2. In a window requiring approximately one third of the time needed to run the orbit ${\bf\Delta}(v,n,m,h)$ remains negative and well separated from zero. Very roughly, this segment starts when the variable $v$ up-crosses the level $-2$ and ends when it up-crosses the level $+5$. On the remaining parts of the orbit, ${\bf\Delta}$ changes sign several times. In particular ${\bf \Delta}$ takes values very close to zero immediately after the top of the spike,  i.e.\ after the variable $v$ has reached its maximum over the stable orbit.

\newpage

\begin{figure}[!h]
\begin{center}
   \includegraphics[width=0.80\textwidth]{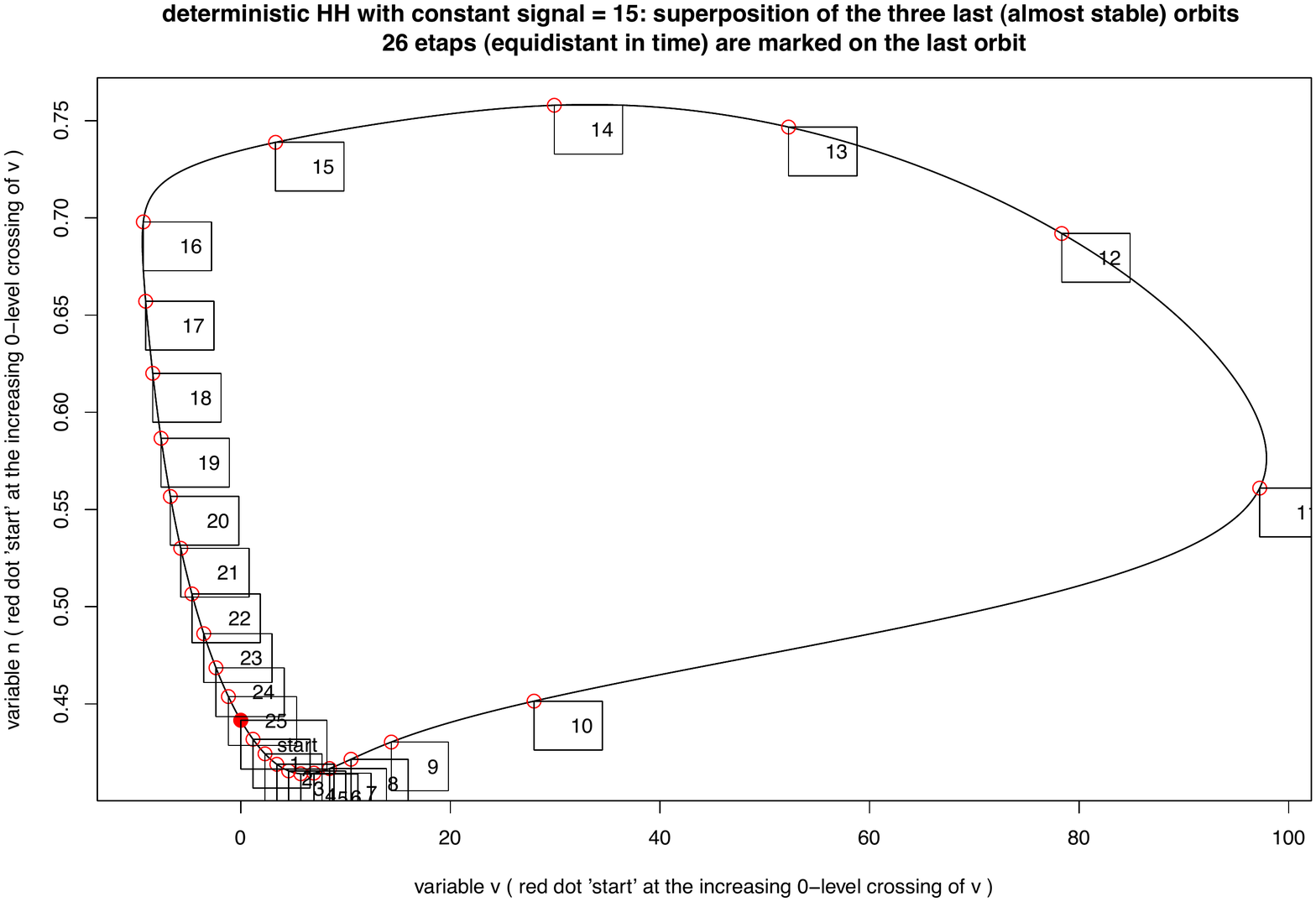} 
\end{center}
\caption{\small Stable orbit of the deterministic system (HH) with constant input $c=15$. }
\end{figure}

\begin{figure}[!h]
\begin{center}
   \includegraphics[width=0.80\textwidth]{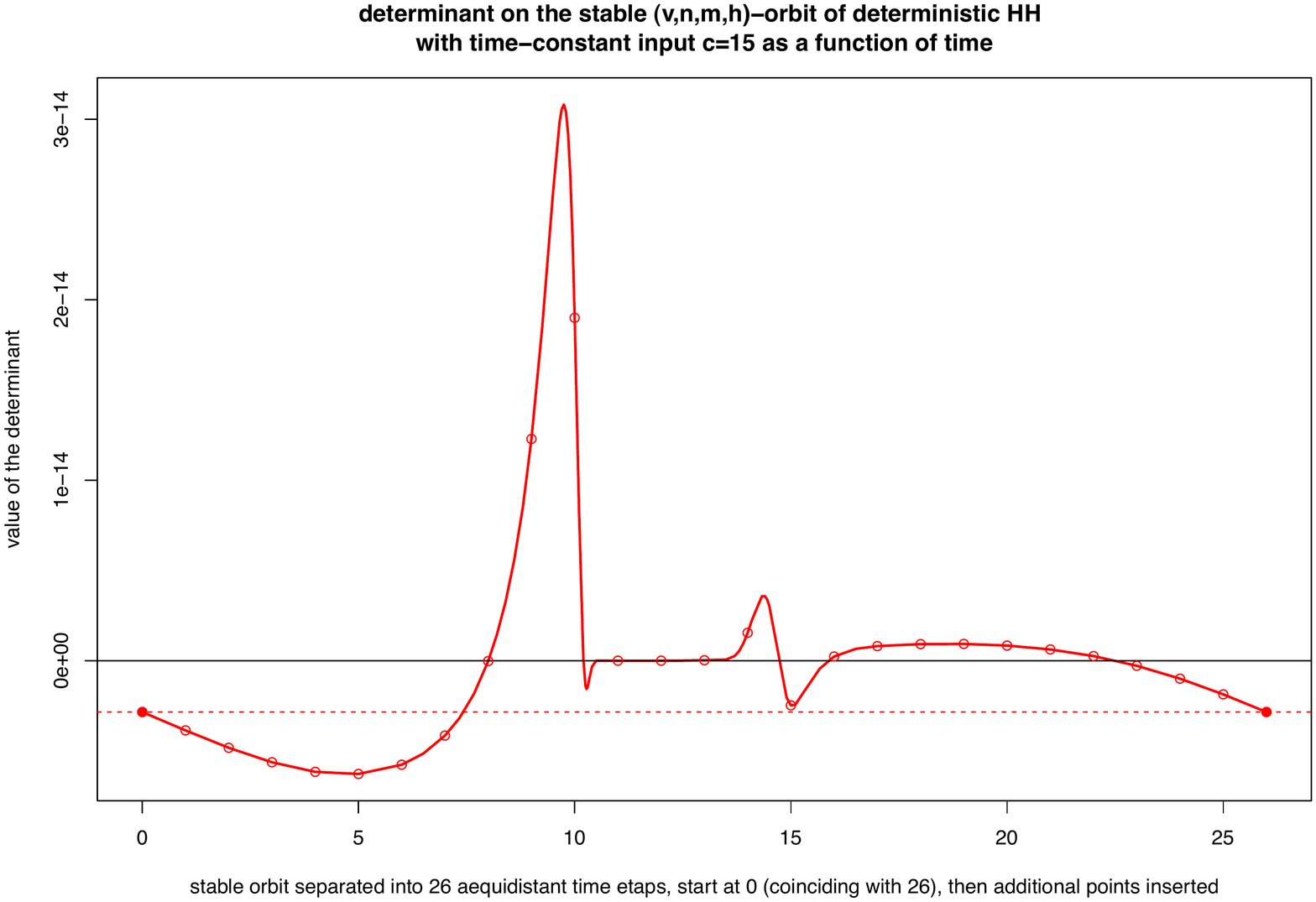} 
\end{center}
\caption{\small Determinant ${\bf \Delta}$ calculated on the stable orbit of the deterministic system (HH) with constant input $c=15$. The time needed to run the orbit is $\approx 12.56$ ms.}
\begin{center}
   \includegraphics[width=0.8\textwidth]{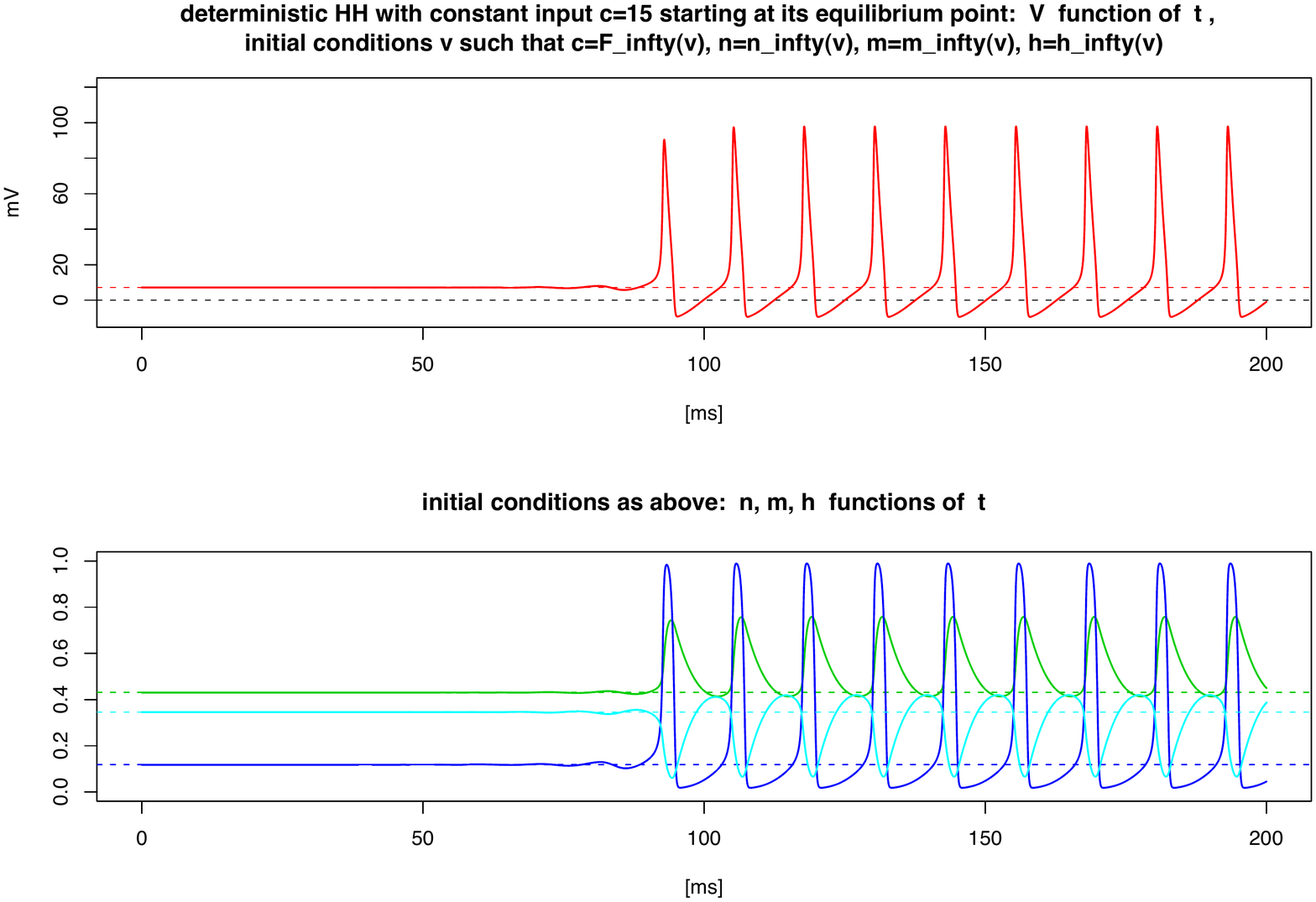} 
\end{center}
\caption{\small Deterministic HH with constant input $c=15.$
}
\end{figure}

\subsection{Positivity regions for (\ref{HHrandominput})}
In this section we apply the results of section \ref{positivedensities} to (\ref{HHrandominput}). Remember that by comparing (\ref{internal1}) and (\ref{HHrandominput}) we see that  $\xi_t-\xi_0$ corresponds to $\int_0^t I(s)ds$.

\noindent Consider first suitable constant $I(t)\equiv c$, fix $\zeta\in U$ and $t>0$, and consider 
\begin{equation}\label{equilibria}
x_c := ( v_c, n_\infty(v_c), m_\infty(v_c), h_\infty(v_c), \zeta ) \quad,\quad x'_c := ( v_c, n_\infty(v_c), m_\infty(v_c), h_\infty(v_c), \zeta+ct )
\end{equation}
where $v_c$ is the unique solution of $F(v_c,n_\infty(v_c), m_\infty(v_c), h_\infty(v_c))=c$ (see section \ref{ex:3}). Let us denote by $P_{s,t}(\cdot,\cdot)_{s<t}$ the semigroup of the process $(X_t)_{ t \geq 0} $ which satisfies (\ref{HHrandominput}). Then Propositions \ref{theo:posinternal} and \ref{cor:3} read as follows.

\begin{prop}\label{cor:1} Assume that $\zeta+cs\in U$ for all $0\le s\le t$. Consider $x_c$ and $x'_c$ defined in (\ref{equilibria}). Then for all $\vep>0$, there exists $\delta > 0$ such that for all $ x'' \in B_\delta ( x_c) , $ $    P_{ 0, t} ( x'' , B_{\varepsilon }( x_c') ) > 0$.
\end{prop} 
\begin{prop}\label{cor:2} 
We keep the assumptions and notations of Proposition \ref{cor:1} and we assume moreover that ${\bf \Delta}( v_c, n_\infty(v_c), m_\infty(v_c), h_\infty(v_c) )\neq 0$. There exists $\delta>0$ such that 
for $ K_c = \overline B_\delta ( x_c)   $ and $ K_c' = \overline B_\delta (x_c') ,$  
$$ \inf_{ x \in K_c } \inf_{x' \in K_c'}  p_{ 0, t} (x, x') > 0. $$
\end{prop}
Remember that the assumption ${\bf \Delta}( v_c, n_\infty(v_c), m_\infty(v_c), h_\infty(v_c) )\neq 0$ ensures that {\bf (LWH)} holds both at $x_c$ and $x'_c$. We have checked numerically in section \ref{ex:3} that this assumption is satisfied for $c\in ]-6.15, 26.61[$. Hence for this range of $c$ these two propositions apply.

\noindent The second situation that we consider is $I(t)  = a \left( 1 + sin( 2\pi \frac{t}{T} ) \right)$, parametrized by $(a,T)$ where $a>0$ is some constant. There are specified subsets $D_1$, $D_2$, $D_3$, $D_4$ in $(0,\infty)\times(0,\infty)$ with the following properties: for $(a,T)$ in $D_1$ (HH) is periodic with small oscillations which cannot be interpreted as spiking. For  $(a,T)$ in $D_2$ the system moves on a $T$-periodic orbit, and the projection $t\mapsto V_t$ resembles the membrane potential of a regularly spiking neuron (single spikes or spike bursts per orbit). For  $(a,T)$ in $D_3$ the system is periodic with period a multiple of $T$. For $(a,T)$ in $D_4$ it behaves irregularly and does not exhibit periodic behavior (see \cite{Endler}, \cite{RM-80} and \cite{AMI-84}). For parameters $(a,T)\in D_2$, consider the points
$$
x:= (0, n^*, m^*, h^*, \zeta ) \quad,\quad x':= (0, n^*, m^*, h^*, \zeta + \int_0^T I(r) dr )
$$ 
such that $(0, n^*, m^*, h^*)$ corresponds to exactly one point on the stable orbit of (\ref{HH}) at which the membrane potential equals $0$. Our numerical example in Fig.\ 1--3 corresponds to parameter values $(a,T)\in D_2$ such that the solution of (\ref{HH}) performs exactly one tour on the stable orbit during $[0,T]$. 
\begin{prop}\label{cor:4}
Assume that $\zeta+\int_0^t I(r) dr \in U$ for all $0\le t\le T$ and that ${\bf \Delta}(x)\neq 0$. There exists
$ \delta > 0 $ such that for $ K  = \overline B_\delta ( x)   $ and $ K' = \overline B_\delta (x') ,$  
$$ \inf_{ y \in K  } \inf_{y' \in K' }  p_{ 0, T} (y, y') > 0 .$$
\end{prop}
Note that the assumption ${\bf \Delta}( 0, n^*, m^*, h^* )\neq 0$ implies that {\bf (LWH)} holds at both points $x$ and $x'$. In our example of Figure 2 its is satisfied at the point on the orbit at which the membrane potential \ul{up-crosses} level $0$.

\section{Appendix}\label{app}
\subsection{Simple properties of (\ref{internal1})}
{\bf Proof of Proposition \ref{linear}.} Given the trajectory of $X_1 $, the variation of constants method yields \begin{equation}\label{eq:exact}
X_{i,t}=X_{i,0}{\rm e}^{-\int_0^t a_i(X_{1,s})ds}+\int_0^t b_i(X_{1,u}){\rm e}^{-\int_u^t a_i(X_{1,r})dr}du.
\end{equation}
However note that (\ref{eq:exact}) does not provide an explicit formula for $X_{i,t}$ since $X_1 $ depends on $X_i$ (the system is fully coupled). Writing $\int_0^t b_i(X_{1,u}){\rm e}^{-\int_u^t a_i(X_{1,r})dr}du=\int_0^t \frac{b_i(X_{1,u})}{a_i(X_{1,u})}a_i(X_{1,u}){\rm e}^{-\int_u^t a_i(X_{1,u})dr}du$, the assumptions on $a_i(\cdot)$ and $b_i(\cdot)$ imply that
\begin{equation}
0 \leq   X_{i,t}\leq X_{i,0}{\rm e}^{-\int_0^t a_i(X_{1,s})ds}+\int_0^t a_i(X_{1,u}){\rm e}^{-\int_u^t a_i(X_{1,r})dr}du.
\end{equation}
By straightforward integration it follows that
\begin{equation*}
0 \leq  X_{i,t}\leq (X_{i,0}+{\rm e}^{\int_0^t a_i(X_{1,r})dr}-1){\rm e}^{-\int_0^t a_i(X_{1,s})ds}
= 1+(X_{i,0}-1){\rm e}^{-\int_0^t a_i(X_{1,s})ds}.
\end{equation*}
The statement follows. \halmos
\subsection{Proof of (\ref{eq:classical2})}
We keep the notations introduced in the proof of Theorem \ref{theo:main2} as well as in section \ref{sec:hoerm}. In order to establish (\ref{eq:classical2}) we extend the argument of \cite{stefano}, Theorem 2.3. To sum up this argument we can say that by an iterative procedure on the Sobolev norms of $  {\cal H} ( \XXX_{t-\delta, t} (y) , \Phi ( \XXX_{t-\delta, t} (y) - y_0 ))$ (in the sense of Malliavin calculus) of different indices, it is proved that estimating these Sobolev norms amounts to estimate the Sobolev norms of $\bar X$ and of the inverse of the Malliavin covariance matrix $ (\Gamma_{ \bar X_t})_{i,j} := < D \bar X_{i,t} , D \bar X_{j,t}>_{L^2 [0, t ]} , \, \, 1 \le i, j \le m,$ where $D$ denotes Malliavin derivative. Since by a classical identity this inverse can be written using the inverse of $det \, \, \Gamma_{\XXX_t}$ and the coefficients of $\Gamma$ itself, the key ingredient is finally to estimate the Sobolev norms of $\bar X$ and expressions of the form $E_z \left( | det \, \, \Gamma_{\XXX_t}|^{-p} \right)^{1/p}$. 
We show below that no difficulty comes from the Sobolev norms of $\bar X$ and we prove that for any $p\geq 1$ and $t\leq 1$, for any $N\in \bbn$ and $z$ such that ${\rm dim}\, \, {\rm LA}(\call_N)(s,z)=m,\, \, \, \forall s \in [0,t]$, 
\begin{equation}\label{eqimportant3}
E_z \left( | det \, \, \Gamma_{\XXX_t}|^{-p} \right)^{1/p} \le C(p,m, N, z) {\, \, t ^{-m(1+N)} }.
\end{equation}
Formula (\ref{eqimportant3}) is the main step to obtain (\ref{eq:classical2}). Indeed it suffices to apply it to the process $\bar X_{t-\delta, t}$ on an interval of length $\delta$ instead of $t$ in (\ref{eqimportant3}). A particular version of (\ref{eqimportant3}) obtained by taking $N=0$ is proved in \cite{stefano} where the restriction to $N=0$ is possible due to the fact that local ellipticity is assumed to hold. However local ellipticity fails to hold in our framework. This is why we prove the general version of (\ref{eqimportant3}). 

\noindent We proceed in three steps. In the first step we check that the usual upper bound for the Sobolev norms of $\bar X$ is still valid and at the end of this step we obtain an expression of a key term of $\Gamma_{ \bar X_t}$ that involves the successive Lie brackets introduced in section \ref{sec:hoerm}. The scheme of this argument is classical (cf. \cite{KS}) but we have to take care of the time dependence in the drift. We describe its main points for the sake of completeness. In the second step we prove (\ref{eqimportant3}) where $N$ is the order of the successive Lie brackets that we need to generate ${\bbr}^m$ according to {\bf (LWH)}. When local ellipticity holds, the diffusion coefficients themselves generate ${\bbr}^m$ and it is not necessary to compute Lie brackets ($N=0$). Finally, in the third step, we show how the arguments of the proof of Theorem 2.3 of \cite{stefano} allow to obtain (\ref{eq:classical2}) from (\ref{eqimportant3}), with $t = \delta$.

\noindent {\it Step 1.} Let 
\begin{equation*}
 \tilde{\bar b}_i (t,x):= \bar b_i (t,x) - \frac12 \sum_{k=1}^m \bar\sigma_k (x) \frac{ \partial \bar \sigma_i }{\partial x_k } (x) ,\;  1 \le i \le m,
 \end{equation*}
be the Stratonovich drift for (\ref{eq:processgood}) and $\bar A_0  := \frac{\partial}{\partial t } + \tilde{\bar b} $, $ \bar A_1 := \bar \sigma$, the corresponding vector fields.
Define $(Y_t)_{i,j} := \frac{ \partial \bar X_{i,t}}{\partial x_j} , 1 \leq i, j \leq m$. Then $Y$ satisfies the following linear SDE with bounded coefficients w.r.t. time and space,
$$ Y_t = I_m + \int_0^t  \partial \bar b ( s, \bar X_s) Y_s  ds + \int _0^t \partial \bar \sigma ( \bar X_s) Y_s d W_s ,$$
where $I_m$ is the $m\times m-$unity matrix and $ \partial \bar b $ and $ \partial \bar \sigma $ are the $m\times m-$matrices having components $ ( \partial \bar b)_{i,j} (t, x) = \frac{\partial \bar b_i }{\partial x_j} (t, x ) $ and  $ ( \partial \bar \sigma)_{i,j} ( x) = \frac{\partial \bar \sigma_i }{\partial x_j} ( x )  .$
By means of It\^o's formula, one shows that $Y_t$ is invertible. Its inverse $Z_t$ satisfies the linear SDE (again with bounded coefficients w.r.t. time and space) given by 
\begin{equation}\label{eq:z}
Z_t = I_m - \int_0^t \partial \tilde{\bar b} (s,\bar X_s)  Z_s ds - \int_0^t  \partial \bar \sigma ( \bar X_s) Z_s \circ dW_s   ,
\end{equation}
where $ \circ d W_s $ denotes the Stratonovich integral. In this framework, the following estimates are classical (see e.g. \cite{KS}) and will be sufficient for our purpose. For all $ 0 \le s \le t \le T ,$ for all $ p \geq 1,$
\begin{equation}\label{eq:ub2}
\sup_{ s \le t } E ( |(Z_s)_{i,j}|^p) \le C(T,p,m, \bar b , \bar \sigma ) ,\; 1 \le i, j \le m ,
\end{equation} 
and 
\begin{equation}\label{eq:ub3}
\sup_{ r_1 , \ldots , r_k \le t } E \left(  | D_{r_1, \ldots , r_k} \bar X_{i,t}|^p  \right) \le C(T,p,m,k, \bar b , \bar \sigma )  \left( t^{1/2}  +  1\right)^{ (k+1)^2 p} ,
\end{equation}
where the constants $C(T,p,m,k, \bar b , \bar \sigma )$ depend only on the bounds of the space derivatives of $\bar b$ and $\bar \sigma.$ Up to this point, the fact that the drift coefficient depends on time did not play an important role since all coefficients are bounded, uniformly in time. 

\noindent {\it Step 2.} It is well known (see for example \cite{Nualart}, page 110, formula (240)) that 
$$ \Gamma_{ \bar X_t} = Y_t \left( \int_0^t Z_s \bar \sigma ( \bar X_s) \bar \sigma^* ( \bar X_s) Z_s^* ds \right) Y_t^* .$$
In order to prove (\ref{eqimportant3}) one has to evaluate the latter integral and therefore to control expressions of the form 
$ Z_s V(s, \XXX_s) ,$ where $V (t,x) $ is a smooth function. This is done by iterating the formula,
\begin{eqnarray}\label{eq:important2}
Z_t V (t, \XXX_t ) &=& V(0, x) + \int_0^t Z_s [ \bar \sigma , V ] ( s, \XXX_s) \circ d W_s 
+ \int_0^t Z_s  [ \frac{\partial}{\partial t} + \tilde{ \bar b} , V] (s, \XXX_s) ds\\
&=& V(0, x) + \int_0^t Z_s [ \bar A_1 , V ] ( s, \XXX_s) \circ d W_s 
+ \int_0^t Z_s  [ \bar A_0 , V] (s, \XXX_s) ds\nonumber, 
\end{eqnarray}
starting with $V\equiv \bar\sigma$, where we identify functions with vector fields (cf. \cite{Nualart}, formula (2.42)). Here, the fact that the drift coefficient is time dependent is important and gives rise to the extra term $ \frac{\partial}{\partial t} $ within the second integral of the first line. In particular with $V\equiv \bar\sigma$ we obtain (cf. (1.9) of \cite{C-M})
\begin{equation*}
Z_t \bar\sigma (t, \XXX_t ) = \bar\sigma(x) + \int_0^t Z_s  [ \bar A_0 , \bar A_1] (s, \XXX_s) ds.\end{equation*}
Iterating (\ref{eq:important2}) we see that $Z_s \bar \sigma (\XXX_s) $ can be written as the sum of two terms. The first term is a finite sum of iterated It\^o integrals where the integrands are $\bar A_1$ and the successive Lie brackets of order at most $N$ obtained with $\bar A_1$ and $\bar A_0$. The second term is a remainder $R_N$ (this is analogous to Theorem 2.12 of \cite{KS}). The most important feature is that the behavior of $R_N$ depends only on the supremum norms of derivatives with respect to time and space of $\bar b $ and with respect to space of $\bar \sigma .$ Based on (\ref{eq:important2}), (\ref{eqimportant3}) follows by Theorem (2.17), estimate (2.18) of \cite{KS}.

\noindent {\it Step 3.} Once (\ref{eqimportant3}) is established, (\ref{eq:classical2}) follows by a straightforward adaptation of the proof of Theorem 2.3 of \cite{stefano}. For completeness let us note that (\ref{eq:ub3}) is the same bound as (2.17) in \cite{stefano} 
whereas (\ref{eqimportant3}) plays the role of (2.20) in \cite{stefano}. For $t$ close to zero, the right-hand side of (2.20) in \cite{stefano} is of order $t^{-m}$ due to the local ellipticity condition, while our bound is of order  $t^{-m(1+N)}$ due to our condition {\bf (LWH)}. Plugging (\ref{eqimportant3}) and (\ref{eq:ub3}) in (2.25) of \cite{stefano} (cf.\ the proof of (2.23)) replaces the r.h.s. obtained there by $O\left( t^{-p\left\{ mN + 1 \right\}} \right)$ for small $t$. 
With such changes, the argument developed there goes through. In our framework we end up with $ O\left( t^{-m\, k_N} \right) $ for small $t$ as r.h.s., with some positive constant $k_N$ depending on {\bf (LWH}). \halmos


\section*{Acknowledgments}
We thank Vlad Bally and Michel Bena\"im for very stimulating discussions.

\end{document}